\documentclass[12pt,amstex]{article}
\usepackage{color}

\usepackage{latexsym}
\usepackage{amsthm}

\usepackage{amsmath}
\usepackage{amsfonts}
\usepackage{amssymb}

\input xy
\xyoption{all}

\usepackage{graphicx}
\usepackage[curve]{xypic}
\newtheorem{thm}{Theorem}[section]
\newtheorem{sub}[thm]{}
\numberwithin{equation}{section}
\newcommand{\cblue}{\color{black}}
\newcommand{\cbb}{\color{black}}
\newcommand{\gl}{\lambda}
\newcommand{\ga}{\alpha}
\newcommand{\gb}{\beta}
\newcommand{\sq}{\sqsubseteq}
\newcommand{\h}{\hspace*{-\parindent}}

\newcommand{\ra}{\rightarrow}

\newcommand{\A}{\mathcal{A}}

\newcommand{\Ij}{Injectivity }

\newcommand{\cat}{category }

\newcommand{\mor}{morphism}

\newcommand{\ijw}{injective w.r.t. }
\newcommand{\ps}{pushout }
\newcommand{\bl}{\color{black}}
\newcommand{\bb}{\color{black}}
\newcommand{\bs}{\begin{sub} }
\newcommand{\es}{\end{sub} }
\newcommand{\ha}{\hat{A}}
\newcommand{\ovh}{\overline{\mathcal{H}}}
\newcommand{\wth}{\widehat{\mathcal{H}}}

\newcommand\id{\operatorname{id}}

\newcommand\cd{\mathcal {D}}
\newcommand\ce{\mathcal {E}}
\newcommand\ca{\mathcal {A}}
\newcommand\cb{\mathcal {B}}

\newcommand\ch{\mathcal {H}}

\newcommand\cm{\mathcal {M}}

\newcommand\Set{\mathbf{Set}}

\newcommand\theo{\textbf{Theorem}}
\newcommand\cor{\textbf{Corollary}}
\newcommand\lem{\textbf{Lemma}}
\newcommand\prop{\textbf{Proposition}}
\newcommand\defn{\textbf{Definition}}
\newcommand\rem{\textbf{Remark}}

\newcommand\exa{\textbf{Example}}
\newcommand\exas{\textbf{Examples}}
\newcommand\pf{\vspace*{2mm}

\hspace*{-\parindent}\textbf{Proof}}
\newcommand\m{\prime}
\newcommand\ova{\overline{A}}
\input xy
\xyoption{all}

\long\def\symbolfootnote[#1]#2{\begingroup%
\def\thefootnote{\fnsymbol{footnote}}\footnote[#1]{#2}\endgroup}

\begin{document}
\title{A Logic of Injectivity}  
\author{J.  Ad\'amek,
M. H\'ebert  and L. Sousa\footnote{The third author acknowledges
financial
    support by  the Center of Mathematics
  of the University of Coimbra and the School of Technology of Viseu}
}


\date{}



\maketitle

\begin{abstract}
Injectivity of objects with respect to a set $\ch$ of morphisms is
an important concept of algebra, model theory and homotopy theory.
Here we study the logic of injectivity consequences of $\ch$, by
which we understand morphisms $h$ such that injectivity with respect
to $\ch$ implies injectivity with respect to $h$. We formulate three
simple deduction rules for the injectivity logic and for its
finitary version where \mor s between finitely ranked objects are
considered only, and prove that they are sound in all categories,
and complete in all ``reasonable" categories.
\end{abstract}


 \section{Introduction}

 Recall that an object $A$ is \ijw a morphism
$h:P\ra P^\m$ provided that every morphism  from $P$ to $A$ factors
through $h$. We address the following problem: given a set
$\mathcal{H}$ of morphisms, which morphisms $h$ are
\textit{injectivity consequences} of $\mathcal{H}$ in the sense that
every object \ijw all members of $\mathcal{H}$ is also \ijw $h$? We
denote the injectivity consequence relationship by
$\mathcal{H}\models h$.

This is a classical topic in general algebra: the
\textit{equational logic} of Garrett Birkhoff \cite{B} is a
special case. In fact, an equation $s=t$ is a pair of elements of
a free algebra $F$, and that pair generates a congruence $\sim$ on
$F$. An algebra $A$ satisfies $s=t$ iff it is \ijw the canonical
epimorphism $$h: F\ra F/\sim.$$ Thus, if we restrict our sets
$\mathcal{H}$ to regular epimorphisms with free domains, then the
logic of injectivity becomes precisely the equational logic.
However, there are other important cases in algebra: recall for
example the concept of injective module, where $\mathcal{H}$ is
the set of all monomorphisms (in the category of modules).

To mention an example from homotopy theory, recall that a
\textit{Kan complex} \cite{K} is a simplicial set injective w.r.t.
all the monomorphisms $\Delta^k_n\hookrightarrow \Delta_n$ (for
$n,\, k \in \mathbb{N}, \, k\leq n$) where $\Delta_n$ is the complex
generated by a single $n$-simplex  and $\Delta^k_n$ is the
subcomplex obtained by deleting the $k$-th 1-simplex  and all
adjacent faces. We can ask for example whether Kan complexes can be
specified by a simpler collection of monomorphisms, as a special
case of our injectivity logic.

Injectivity establishes a Galois correspondence between objects and
morphisms of a category. The closed families on the side of objects
are called \textit{injectivity classes}: for every set $\ch$ of \mor
s we obtain the injectivity class Inj$\ch$, i.e., the class of all
objects \ijw $\ch$. In \cite{AR2} small-injectivity classes in
locally presentable categories were characterized as precisely the
full accessible subcategories closed under products, and in
\cite{RAB} this was sharpened in the following sense. Let us call a
morphism \textit{$\lambda$-ary} if its domain and codomain are
$\lambda$-presentable objects. Injectivity classes with respect to
$\lambda$-ary morphisms are precisely the full subcategories closed
under products, $\gl$-filtered colimits, and $\lambda$-pure
subobjects. For injectivity w.r.t. cones or trees of morphisms
similar results are in \cite{AN77} and \cite{NS77}.

In the present paper we study closed sets on the side of morphisms,
i.e., we develop a deduction system for the above  injectivity
consequence relationship $\models$. It has altogether  three
deduction rules, which are quite intuitive. Firstly, observe that
every object injective w.r.t. a composite $h=h_2\cdot h_1$  is
injective w.r.t. the first morphism $h_1$. This gives us the first
deduction rule

\vspace*{4.5mm} \hspace*{1cm}\begin{tabular}{p{2.7cm}l}{\sc
cancellation}\\ \\
\end{tabular}\hspace*{6mm}\begin{tabular}{c}$h_2 \cdot h_1$
\\ \hline $ h_1$ \\ \\ \end{tabular}

\hspace*{-\parindent}It is also easy to see that injectivity
w.r.t. $h$ implies injectivity w.r.t. any morphism $h^\m$ opposite
to $h$ in a pushout (along an arbitrary morphism), which yields
the rule

\hspace*{1cm}\begin{tabular}{p{2.7cm}l}{\sc pushout}\\
\end{tabular}\hspace*{3mm}\begin{tabular}{l}$ h$\\ \hline
 $ h'$ \\ \end{tabular}
\begin{tabular}{l}\\ {\hspace*{10mm}for every pushout }\\ \\ \end{tabular}
\begin{tabular}{l}$
\xy (0,0)*{\xymatrix{\ar[r]^h\ar[d]& \ar[d]\\
\ar[r]^{h'}& }}="D"; (7,-10.5)*{}="A"; (7,-7)*{}="B";
(10.5,-7)*{}="C"; "A"; "B" **\dir{-}; "B"; "C" **\dir{-};
\endxy
$\end{tabular}

\vspace*{3mm}

\hspace*{-\parindent}Finally, an object injective w.r.t. two
composable morphisms is also injective w.r.t. their composite. The
same holds for three, four, $\dots$ morphisms --  but also for a
\textit{transfinite composite} as used in homotopy theory. For
example, given an $\omega$-chain of morphisms
$$\xymatrix{A_0\ar[r]^{h_0}&A_1\ar[r]^{h_1}&A_2\ar[r]^{h_2}&\dots}$$
then their $\omega$-composite is the first morphism $c_0:A_0\ra C$
of (any) colimit cocone $c_n:A_n\ra C\, (n\in \mathbb{N})$ of the
chain. Observe that $c_0$ is indeed an injectivity consequence of
$\{h_i;\, i<\omega \}$. For every ordinal $\lambda$ we have the
concept of a $\lambda$-composite of morphisms (see \ref{a2.9} below)
and the following deduction rule, expressing the fact that an object
injective w.r.t. each $h_i$ is injective w.r.t. the transfinite
composite:

\hspace*{0.6cm}\begin{tabular}{p{5.5cm}l}{\sc transfinite composition}\\
\end{tabular}\hspace*{1.2mm}\begin{tabular}{c}$ h_i\, (i<\lambda)$\\ \hline
 $ h$ \\ \end{tabular}
\begin{tabular}{l}\\ {\hspace*{1.0mm}for every $\lambda$-composite $h$ of $(h_i)_{i<\gl}$ }\\ \\ \end{tabular}

\hspace*{-\parindent}We are going to prove that the Injectivity
Logic based on the above three rules is sound and complete. That
is, given a set $\mathcal{H}$ of morphisms, then
$\mathcal{H}\models h$ holds for precisely those morphisms $h$
which can be proved from assumptions in $\mathcal{H}$ using the
three deduction rules above. This holds in a number of categories,
e.g., in

\begin{enumerate}
\item[(a)] every variety of algebras, \item[(b)] the category of
topological spaces and many nice subcategories (e.g. Hausdorff
spaces), and \item[(c)] every locally presentable category of
Gabriel and Ulmer.
\end{enumerate}
We introduce the concept of a strongly locally ranked category
encompassing (a)-(c) above, and prove the soundness and completeness
of our Injectivity Logic in all such categories.

 Observe that the above logic is
infinitary, in fact, it has a proper class of deduction rules: one
for every ordinal $\lambda$ in the instance of {\sc transfinite
composition}. We also  study, following the footsteps of Grigore
Ro\c su,  the completeness of the corresponding Finitary \Ij Logic:
it is the restriction of the above logic to $\lambda$ finite. Well,
all we need to consider {\cblue are the cases} $\lambda =2$, called
{\sc composition}, and $\lambda=0$, called {\cblue {\sc identity}}:

\vspace*{4.5mm} \hspace*{2cm}\begin{tabular}{p{2.7cm}l}{\sc
composition}\\ \\
\end{tabular}\hspace*{3mm}\begin{tabular}{c}$h_1\; \; h_0$
\\ \hline $ h$ \\ \\ \end{tabular}
\begin{tabular}{l}\\ {\hspace*{10mm}for $h= h_1\cdot h_0$ }\\ \\ \end{tabular}

\hspace*{2cm}\begin{tabular}{p{2.7cm}l}{\cblue {\sc identity}}\\ \\
\end{tabular}\hspace*{3mm}\begin{tabular}{l}\hline $\id_A$\\ \end{tabular}

\h The resulting finitary deductive system (introduced in \cite{ASS}
as a slight modification of the deduction system of Grigore Ro\c su
\cite{R}) has four deduction rules; it is clearly sound, and the
main result of our paper (Theorem \ref{a5.2}) says that it is also
complete with respect to finitary morphisms, i.e., \mor s with
domain and codomain of finite rank. This implies the expected
compactness theorem: every finitary injectivity consequence of a set
$\ch$ of finitary \mor s is an injectivity consequence of some
finite subset of $\ch$.

The completeness theorem for Finitary Injectivity Logic will then
be extended to the $k$-ary Injectivity Logic, defined in the
expected way. Then the full completeness theorem easily follows.

The fact that the full Injectivity Logic above is complete in
strongly locally ranked categories can also be derived from
Quillen's Small Object Argument \cite{Q}, see Remark \ref{b3.9}
below. However our sharpening to the $k$-ary logic for every
cardinal $k$ cannot be derived from that paper, and we consider this
to be a major {\cbb step.}

\vspace*{2mm} \h {\textbf{Related work}} Bernhard Banaschewski and
Horst Herrlich showed thirty years ago that implications in general
algebra can be expressed categorically via injectivity w.r.t.
regular epimorphisms, see \cite{BH}. A generalization to injectivity
w.r.t. cones or even trees of morphisms was studied by Hajnal Andr\'
eka, Istv\'an N\'emeti and Ildik\'o Sain, see e.g. \cite{AN77, AN79,
NS77}.

To see more precisely how that work relates to ours and to classical
logic, consider injectivity in the \cat of all $\Sigma$-structures
(and $\Sigma$-homomorphisms), where $\Sigma$ is any signature. Then
recall from \cite{AR}, 5.33 that there is a natural way to associate
to a (finitary) morphism $f:A\ra B$ a (finitary) sentence
$$f^\m:=\forall X(\wedge A^\m(X)\ra \exists Y (\wedge B^\m(X,Y)))$$
(where $A^\m(X)$ and $B^\m(X,Y)$ are sets of atomic formulas) such
that an object $C$ satisfies $f^\m$ if and only if it is injective
with respect to $f$ (see 2.22 below for more on this).
 Such sentences are called \textit{regular
sentences}.  In this paper we concentrate on the proof theory for
the  (finite and infinite) regular logics. As mentioned above, the
restriction to epimorphisms correspond to considering only the
quasi-equations (i.e., no existential quantifiers), and just
equations if we impose they have projective domains.

Recently, Grigore Ro\c su introduced a deduction system for
injectivity, see \cite{R}, and he proved that the resulting logic is
sound and complete for epimorphisms which are finitely presentable,
see \ref{b3.4}, and have projective domains. A slight modification
of Ro\c su's system was introduced in \cite{ASS}: this is the
deduction system 2.4 below. It differs from \cite{R} by formulating
{\sc pushout}  more generally and using {\sc composition} in place
of Ro\c su's {\sc union}. In \cite{ASS} completeness is proved for
sets of epimorphisms with finitely presentable domains and
codomains. (This is slightly stronger than requiring the
epimorphisms to be finitely presentable, however, without the too
restrictive assumption of projectivity of the domains the logic
fails to be complete for finitely presentable epimorphisms in
general, see \cite{ASS}.)

In the present paper completeness of the finitary logic is proved
for arbitrary morphisms (not necessarily epimorphisms) with finitely
presentable domains and codomains. The fact that the assumption of
epimorphism is dropped makes the proof substantially more difficult.
We present a short proof in locally presentable categories first,
and then a proof of a more general result for strongly locally
ranked categories. We also formulate the appropriate infinitary
logic dealing with arbitrary morphisms.

There are other generalizations of Birkhoff's equational logic which
are, except for the common motivation, not related to our approach.
For example  the categorical approach to logic of (ordered)
many-sorted algebras of Razvan Diaconescu \cite{Dia}, and the logic
of implications in general algebra of Robert Quackenbush \cite{Qua}.

In our joint paper \cite{AHS2} we are taking another route to
generalize the equational logic: we consider orthogonality of
objects to a morphism instead of injectivity. The deduction system
is similar: the rule {\sc cancellation} has to be weakened, and an
additional rule concerning coequalizers is added. We prove the
completeness of the resulting  logic of orthogonality in locally
presentable categories. The corresponding sentences are the so
called limit sentences, $\forall X(\wedge A^\m(X)\ra
\exists!Y(\wedge B^\m(X,Y)))$, where $\exists!Y$ means ``there
exists exactly one $Y$ such that".

\section{Logic of injectivity}

 \setcounter{thm}{-1}

\begin{sub}\label{2.0}{\em {\textbf{Assumption}} Throughout the paper we assume
that we are working in a cocomplete category.}
\end{sub}

\begin{sub}\label{2.1}{\em {\defn} A morphism $h$ is called an
{\cblue  \textit{injectivity consequence}} of a set of morphisms
$\ch$, notation
$$\ch \models h$$
 provided that every object injective w.r.t. all \mor s in $\ch$  is
 also injective w.r.t. $h$.
}
\end{sub}

\begin{sub}\label{2.2}{\em {\exas} (1) A composite $h=h_2\cdot h_1$ is an injectivity consequence of $\{h_1,\, h_2\}$.

(2) Conversely, in every composite $h=h_2\cdot h_1$
 the morphism $h_1$ is an injectivity consequence of $h$: 
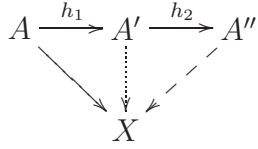
$$\xymatrix{A\ar[r]^{h_1}\ar[dr]&A^\m \ar[r]^{h_2} \ar@{.>}[d]
&A^{\m\m}\ar@{-->}[dl]\\&X&}$$

 (3) In every pushout
 $$
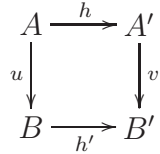
\xymatrix{A\ar[r]^h\ar[d]_u&A^\m\ar[d]^v\\
 B\ar[r]_{h^\m}&B^\m}$$
 $h^\m$ is an injectivity consequence of $h$:
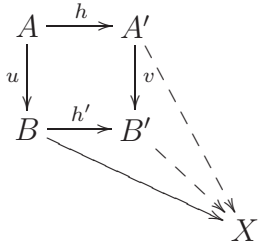
 $$\xymatrix{A\ar[r]^h\ar[d]_u&A^\m\ar[d]^v\ar@{-->}[ddr]&\\
 B\ar[r]^{h^\m}\ar[drr]&B^\m\ar@{-->}[dr]&\\&&X}$$
 }
\end{sub}

\begin{sub}\label{2.3}{\em {\rem} The above examples are exhaustive. More precisely,  the following
deduction system, introduced in \cite{ASS}, see also \cite{R},
(where, however, it was only applied to epimorphisms) will be proved
complete below:}\end{sub}

\begin{sub}\label{2.4}{\em {\defn} The \textit{Finitary Injectivity  Deduction System} consists of one
axiom

\vspace*{6mm}
\hspace*{2cm}\begin{tabular}{p{2.7cm}l}{\sc identity}\\ \\
\end{tabular}\hspace*{3mm}\begin{tabular}{l}\hline $\id_A$\\ \end{tabular}

\vspace*{2mm} and three deduction rules

\vspace*{4.5mm}

\hspace*{2cm}\begin{tabular}{p{2.7cm}l}{\sc composition}
\\ \end{tabular}\hspace*{3mm}\begin{tabular}{ll} $h\; \;
h'$& \\ \hline
 $h' \cdot h$\\ \end{tabular}
 \begin{tabular}{l}{
 if $h^\m\cdot h$ is defined}\end{tabular}

\vspace*{4.5mm} \hspace*{2cm}\begin{tabular}{p{2.7cm}l}{\sc
cancellation}\\ \\
\end{tabular}\hspace*{3mm}\begin{tabular}{c}$h' \cdot h$
\\ \hline $ h$ \\ \\ \end{tabular}

and

\hspace*{2cm}\begin{tabular}{p{2.7cm}l}{\sc pushout}\\
\end{tabular}\hspace*{3mm}\begin{tabular}{l}$ h$\\ \hline
 $ h'$ \\ \end{tabular}
\begin{tabular}{l}\\ {\hspace*{10mm}if }\\ \\ \end{tabular}
\begin{tabular}{l}$
\xy (0,0)*{\xymatrix{\ar[r]^h\ar[d]& \ar[d]\\
\ar[r]_{h'}& }}="D"; (7,-10.5)*{}="A"; (7,-7)*{}="B";
(10.5,-7)*{}="C"; "A"; "B" **\dir{-}; "B"; "C" **\dir{-};
\endxy
$\end{tabular}

\vspace*{6mm}

We say that a morphism $h$ is a \textit{formal consequence} of a set
$\ch$ of morphisms (notation $\ch \vdash h$) in the Finitary
Injectivity Logic if there exists a proof of $h$ from $\ch$ (which
means
 a finite sequence $h_1, \, ..., \, h_n = h$ of morphisms such
that for every $i = 1, ..., n$ the morphism $h_i$ lies in $\ch$ or
is a conclusion of one of the deduction rules whose premises lie in
$\{h_1,...,h_{i-1}\}$).}\es

\bs \label{2.61/2} {\em {\lem}} The Finitary Injectivity Logic is
sound, i.e., if a {\cbb {\mor} } $h$ is a formal consequence of a
set of {\cbb \mor s } $\ch$, then $h$ is an injectivity consequence
of $\ch$. Briefly: $\ch \vdash h$ implies $\ch \models h$.

\vspace*{3mm}

{\em The proof follows from \ref{2.2}.}\es

\bs\label{novo2.7}{\em {\rem} Later  we define finitary morphisms
(as morphisms whose domains and codomains are finitely presentable
(Section 3) or of finite rank (Section 5)), and in Section 6 we
prove that the resulting Finitary Injectivity Logic is complete,
i.e., that
$$\ch \models h\; \mbox{ implies }\; \ch \vdash h$$ for every set
$\ch$ of finitary \mor s and every $h$ finitary.  }\end{sub}

\begin{sub}\label{2.6}{\em {\exa} The following rule

\vspace*{4.5mm}
\hspace*{2cm}\begin{tabular}{p{3.7cm}l}{\sc finite coproduct}\\ \\
\end{tabular}\hspace*{3mm}\begin{tabular}{l}$h_1 \; \; \; h_2 $
\\ \hline $ h_1+h_2$ \\ \\ \end{tabular}\newline
(where for $h_i:A_i\ra B_i$ the morphism $h_1+h_2:A_1+A_2 \ra
B_1+B_2$ is the canonical coproduct morphism) is obviously sound.
Here is a proof in the Finitary Injectivity Logic:

\h Using the pushouts

$$
\xy (0,0)*{\xymatrix{A_1\ar[rr]^{h_1}\ar[d]_{}&& B_1\ar[d]^{}\\
A_1+A_2\ar[rr]_{h_1+\id_{A_2}}&& B_1+A_2}}="D"; (25,-14)*{}="A";
(25,-8)*{}="B"; (35.5,-8)*{}="C"; "A"; "B" **\dir{-}; "B"; "C"
**\dir{-}; (30,0)*{\xymatrix{A_2\ar[rr]^{h_2}\ar[d]_{}&& B_2 \ar[d]^{}\\
B_1+A_2\ar[rr]_{\id_{B_1}+h_2}&&B_1+B_2 }}="D"; (85,-14)*{}="A";
(85,-8)*{}="B"; (95.5,-8)*{}="C"; "A"; "B" **\dir{-}; "B"; "C"
**\dir{-};
\endxy
$$

\vspace*{5mm} \hspace*{-\parindent}we can write

{ \begin{center}\begin{tabular}{p{4.6cm}}
\begin{tabular}{c}$\quad \; \; h_1 \;\; \; \qquad \qquad h_2$\end{tabular}

 \begin{tabular}{p{4cm}}\hline \\ \end{tabular}

\vspace*{-4mm}

  \begin{tabular}{c}$h_1+\id_{A_2} \qquad \id_{B_1}+h_2$\end{tabular}

  \begin{tabular}{p{4cm}}\hline \\ \end{tabular}

\vspace*{-4mm}

\begin{tabular}{c}$\qquad \; \; \; \, h_1+h_2$
\end{tabular}
\end{tabular}
\begin{tabular}{p{4cm}} \vspace*{-5mm} via {\sc pushout}\\
via {\sc composition}
\end{tabular}

\end{center}

\hspace*{-\parindent} since
$h_1+h_2=(\id_{B_1}+h_2)\cdot(h_1+\id_{A_2})$. }}\end{sub}

\begin{sub}\label{2.7}{\em {\exa}   The following rule

\hspace*{1cm}\begin{tabular}{p{4.2cm}l}{\sc finite wide pushout}
\\ \end{tabular}\hspace*{3mm}\begin{tabular}{c} $h_1\; \dots  \;
h_n$ \\ \hline
 $h$\\ \end{tabular}

\h for every wide pushout

\vspace*{-10mm}
$$\xymatrix{&\ar[dl]_{h_1}\ar[d]_{h_2}\ar[dr]_{\dots}^{h_n}&\\
\ar[dr]_{k_1}&\ar[d]_{k_2}&\ar[dl]_{\dots}^{k_n}\\
&C&} \; \; \qquad \mbox{\begin{tabular}{l}\\ \\ \\ \\ where
$h=k_i\cdot h_i$\end{tabular}}$$ is sound. Here is a proof in the
Finitary Injectivity Logic:

If $n=2$ we have

{ \begin{center}
\begin{tabular}{p{2.5cm}}\begin{tabular}{c}$\; \, h_1\qquad h_2$
\end{tabular}

  \begin{tabular}{p{2cm}}\hline \\ \end{tabular}

\vspace*{-4mm}

\begin{tabular}{c}$\qquad k_2$
\end{tabular}

  \begin{tabular}{p{2cm}}\hline \\ \end{tabular}

\vspace*{-4mm}\begin{tabular}{c} $h=k_2\cdot h_2$
\end{tabular} \end{tabular}\begin{tabular}{p{3.8cm}} \vspace*{-5mm}via {\sc pushout}\\
via {\sc composition}
\end{tabular}

\end{center}  }
If $n=3$ denote by $r$ a pushout of $h_1$, $h_2$, then a pushout,
$h_3^\m$,

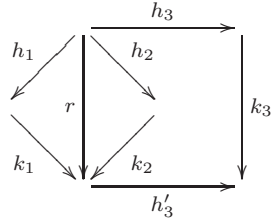
$$\xymatrix{&\ar[dl]_{h_1}\ar[dd]_r\ar[dr]^{h_2}\ar[rr]^{h_3}&&\ar[dd]^{k_3}&\\
\ar[dr]_{k_1}&&\ar[dl]^{k_2}&&\\
&\ar[rr]_{h^\m_3}&&&\\ &&&&}$$

\h of $h_3$ along $r$ forms a wide pushout of $h_1$, $h_2$ and
$h_3$:

{ \begin{center} \begin{tabular}{p{2.5cm}}
\begin{tabular}{c}$\, h_1\; \;
\; h_2\; \; \; h_3$\end{tabular}

  \begin{tabular}{p{2cm}}\hline \\ \end{tabular}

\vspace*{-4mm} \begin{tabular}{c}$\qquad k_2$ \end{tabular}

  \begin{tabular}{p{2cm}}\hline \\ \end{tabular}

\vspace*{-4mm} \begin{tabular}{c}$\qquad r$\end{tabular}

  \begin{tabular}{p{2cm}}\hline \\ \end{tabular}

\vspace*{-4mm} \begin{tabular}{c}$\qquad k_3$ \end{tabular}

  \begin{tabular}{p{2cm}}\hline \\ \end{tabular}

\vspace*{-4mm} \begin{tabular}{c} $\; h=k_3 \cdot h_3$
\end{tabular} \end{tabular} \begin{tabular}{p{4cm}} \vspace*{-5mm}\begin{tabular}{l} via {\sc
pushout}\end{tabular}

\vspace*{2.1mm}

\begin{tabular}{l} via {\sc composition}\end{tabular}

\vspace*{2.1mm}

  \begin{tabular}{l} via {\sc
 pushout}\end{tabular}

 \vspace*{2.1mm}

  \begin{tabular}{l}via {\sc composition}
\end{tabular}
\end{tabular}
\end{center}  }

\h Etc.

}\end{sub}

\begin{sub}\label{a2.8}{\em {\rem} We want to define a composition of
a chain of $\lambda$ morphisms for every ordinal
 $\gl$ (see the case
$\gl=\omega$ in the Introduction). Recall that a
\textit{$\gl$-chain} is a functor $A$ from $\gl$, the well-ordered
category of all ordinals $i<\gl$.

Recall further that $\gl^+$ denotes the successor ordinal, i.e., the
set of all $i\leq \gl$.}\end{sub}

\begin{sub}\label{a2.9}{\em {\defn} (i) We call a $\gl$-chain $A$  \textit{smooth} if
for every limit ordinal $i<\gl$  we have
$$A_i=\mbox{co}\hspace*{-0.8mm}\lim_{j<i} A_j$$ with the colimit
cocone of all $a_{ji}=A(j\ra i)$.

(ii) A morphism $h$ is called a \textit{$\gl$-composite} of
morphisms $(h_i)_{i<\gl}$, where $\gl$ is an ordinal, if there
exists a smooth $\gl^+$-chain $A$ with connecting \mor s
$a_{ij}:A_i\ra A_j$ for $i\leq j\leq \gl$ such that
$$h_i=a_{i,i+1}\; \; \; \mbox{for all $i<\gl$}$$
and $$h=a_{0,\gl}.$$  }\end{sub}

\begin{sub}\label{a2.10}{\em {\exas} $\gl=0$: No morphism $h_i$ is given,
just an object $A_0$; and $h=a_{0,0}$ is the identity morphism of
$A_0$.

$\gl=1$: A morphism $h_0$ is given, and  we have $h=a_{0,1}=h_0$.
Thus, a 1-composite of $h_0$ is $h_0$.

$\gl=2$: This is the usual concept of composition: given morphisms
$h_0$, $h_1$, their 2-composite exists iff they are composable.
Then $h_1\cdot h_0$ is the 2-composite.

$\gl=\omega$: This is the case mentioned in the Introduction.
Observe that, unlike the previous cases, an $\omega$-composite is
only unique up to  isomorphism. }\end{sub}

\begin{sub}\label{aa2.10} {\em \lem} A $\gl$-composite of
\mor s $(h_i)_{i<\gl}$ is an injectivity consequence of these
morphisms.

{\em {\pf} This is a trivial transfinite induction on $\gl$. In case
$\gl=0$ this states that $\id_A$ is an injectivity consequence of
$\emptyset$, etc.}
\end{sub}

\begin{sub}\label{aa2.11}{\em {\defn} The \textit{Injectivity Deduction
System} consists of the deduction rules

\vspace*{3mm}

\hspace*{0.6cm}\begin{tabular}{p{2.7cm}l}{\sc cancellation}\\ \\
\end{tabular}\hspace*{3mm}\begin{tabular}{c}$h' \cdot h$
\\ \hline $ h$ \\ \\ \end{tabular}

\hspace*{0.6cm}\begin{tabular}{p{2.7cm}l}{\sc pushout}\\
\end{tabular}\hspace*{3mm}\begin{tabular}{c}$ h$\\ \hline
 $ h'$ \\ \end{tabular}
\begin{tabular}{l}\\ {\hspace*{10mm}for every pushout }\\ \\ \end{tabular}
\begin{tabular}{l}$
\xy (0,0)*{\xymatrix{\ar[r]^h\ar[d]& \ar[d]\\
\ar[r]^{h'}& }}="D"; (7,-10.5)*{}="A"; (7,-7)*{}="B";
(10.5,-7)*{}="C"; "A"; "B" **\dir{-}; "B"; "C" **\dir{-};
\endxy
$\end{tabular}

\vspace*{3mm}

\h and the rule scheme (one rule for every ordinal $\gl$)

\hspace*{0.6cm}\begin{tabular}{p{5.2cm}l}{\sc transfinite composition}\\
\end{tabular}\hspace*{3mm}\begin{tabular}{c}$ h_i \; (i<\gl)$\\ \hline
 $ h$ \\ \end{tabular}
\begin{tabular}{l}\\ {for every $\gl$-composite $h$ of $(h_i)_{i<\gl}$ }\\ \\ \end{tabular}

We say that a morphism $h$ is a \textit{formal consequence} of a set
$\ch$ of morphisms (notation $\ch \vdash h$) in the Injectivity
Logic if there exists a proof of $h$ from $\ch$ (which means  a
chain $(h_i)_{i \leq n}$ of morphisms, where $n$ is an ordinal, such
that $h = h_n$, and each $h_i$ either lies in $\ch$, or is a
conclusion of one of the deduction rules whose premises lie in
$\{h_j\}_{j < i}$). }\end{sub}

\begin{sub}\label{aa2.13}{\em {\lem} } The \Ij
Logic is sound, i.e., if a {\mor} $h$ is a formal consequence of a
set $\ch$ of \mor s, then $h$ is an injectivity consequence of
$\ch$. Briefly: $\ch \vdash h \; \mbox{ implies }\; \ch \models h.$

\vspace*{2.5mm}

{\em The proof (using \ref{aa2.10}) is elementary.}\end{sub}

\begin{sub}\label{aa2.14}{\em {\rem}  In \ref{aa2.11} we can replace
 {\sc transfinite composition} by the deduction rule
{\sc wide pushout}, see below, which makes use of the (obvious) fact
that an object $A$ \ijw a set $\{h_i\}_{i<\gl}$ of morphisms having
a common domain is also \ijw their wide pushout. Let us note here
that this rule does not replace {\sc pushout} of \ref{aa2.11}
(because in the latter a pushout of $h$ along an {\it arbitrary}
{\mor} is considered). }\end{sub}

\begin{sub}\label{aa2.15}{\em {\defn}  The deduction rule

\vspace*{3mm}

\hspace*{0.6cm}\begin{tabular}{p{3.7cm}l}{\sc wide pushout}\\ \\
\end{tabular}\hspace*{3mm}\begin{tabular}{c}$h_i\; (i<\gl) $
\\ \hline $ h$ \\  \end{tabular}\begin{tabular}{l}\hspace*{3mm}for $h$ a wide pushout of $\{h_i\}_{i<\gl}$
\end{tabular}

\h applies, for every cardinal $\gl$, to an arbitrary object $P$
and an arbitrary set $\{h_i\}$ of $\gl$ morphisms with the common
domain $P$ and the following wide pushout

\vspace*{-1cm}
$$\xymatrix{&P\ar[ld]_{h_i}\ar[d]\ar[rd]&\\
P_i\ar[rd]_{k_i}&\ar[d]&\mbox{\hspace*{8mm}$\dots$}\ar[ld]\\&Q&&}\begin{array}{l}
   \\
   \\
   \\
   \\
   \\
  h=k_i\cdot
h_i \mbox{ (for any $i$)} \\
\end{array}
 $$

{\rem}  Again, this is a scheme of deduction rules: for every
{\cbb cardinal} $\gl$ we have one rule $\gl$-{\sc wide pushout}.
Observe that $\gl=0$ yields the rule {\sc identity}.}\end{sub}

\begin{sub}\label{aa2.16} {\em {\lem}  } The \Ij Deduction System \ref{aa2.11} is equivalent to the deduction system

\centerline{{\sc composition, cancellation, pushout} and {\sc wide
pushout}.}

{\em {\pf} (1) We can derive {\sc wide pushout} from \ref{aa2.11}.
For  every ordinal number $\gl$ we derive the rule

\hspace*{3.5cm}\begin{tabular}{c}\\$h_i\; (i<\gl) $
\\ \hline $ h$ \\ \\ \end{tabular}\hspace*{6mm}\begin{tabular}{l}for $h$ a wide pushout of $\{h_i\}_{i<\gl}$
\end{tabular}

\h by transfinite induction on the ordinal $\gl$. We are given an
object $P$ and \mor s $h_i:P\ra P_i\, (i<\gl)$. The case $\gl=0$
is trivial, from $\gl$ derive $\gl+1$ by using {\sc pushout}, and
for limit ordinals $\gl$ form the restricted multiple pushouts
$Q_j$ of \mor s $h_i$ for $i<j$, and observe that they form a
smooth chain whose composite is a multiple pushout of all $h_i$'s.

(2) From the system in \ref{aa2.16} we can derive the rule
$\gl$-{\sc composition}, where $\gl$ is an arbitrary ordinal: the
case $\gl=0$ follows from  0-{\sc wide pushout}. The isolated step
uses {\sc composition}: the $(\gl+1)$-composite of $(h_i)_{i\leq
\gl}$ is simply $h_\gl \cdot k$ where $k$ is the $\gl$-composite of
$(h_i)_{i < \gl}$. In the limit case, use the fact that a composite
$h$ of $(h_i)_{i < \gl}$ is a wide pushout of $\{k_i\}_ {i<\gl}$,
where $k_i$ is a composite of $(h_j)_{j < i}$.}\end{sub}

\begin{sub}\label{ver} {\em {\rem } For every infinite cardinal $k$ the \textit{$k$-ary Injectivity Deduction
System} is the system \ref{aa2.11} where $\lambda$ ranges through
ordinals smaller than $k$. A proof of  a morphism $h$ from a set
$\ch$ in the $k$-ary Injectivity Logic is, then, a proof of length
$n < k$ using only the deduction rules with $\gl$ restricted  as
above.  The last lemma can, obviously, be formulated under this
restriction in case we use the scheme
  $\gl$-{\sc wide pushout}
  for all cardinals $\gl<k$.}\end{sub}

\begin{sub}\label{aa2.17}{\em {\defn} The deduction rule

\vspace*{5mm}

 \hspace*{1cm}\begin{tabular}{p{4.2cm}l}{\sc coproduct}
\\ \end{tabular}\hspace*{3mm}\begin{tabular}{c} $h_i\; (i<\gl)$ \\
$\xy (0,0)*{}="A"; (20,0)*{}="B"; "A"; "B" **\dir{-} \endxy$\\
 $
 {\coprod_{i<\gl} h_i}$\\ \end{tabular}

\vspace*{3mm}

\h applies, for every cardinal $\gl$, to an arbitrary collection
of $\gl$ \mor s $h_i:A_i\ra B_i$.

}\end{sub}

\begin{sub}\label{aa2.18}{\em {\lem}} The Injectivity Deduction
System \ref{aa2.11} is equivalent to the deduction system of
\ref{aa2.16}
 with {\sc wide pushout} replaced by

 \centerline{\sc {\cblue identity}  \hspace*{1mm} $+$ \hspace*{1mm} coproduct}
{\em {\pf} (1) {\sc coproduct} follows from \ref{aa2.16}. In fact,
${\coprod_{i<\gl}h_i:\coprod_{i<\gl}A_i\ra \coprod_{i<\gl}B_i}\,$
is a wide pushout of the \mor s $ {k_j:\coprod_{i<\gl}A_i\ra
\coprod_{i< j }A_i +B_j+ \coprod_{j< i<\gl}A_i}$, where $j$ ranges
through $\gl$, with components $\, \id_{A_i}\, (i\not= j)$ and
$h_j$, and $k_j$ is a pushout of $h_j$ along the $j$-th coproduct
injection of $ {\coprod_{i<\gl}A_i}$.

(2) Conversely, {\sc wide pushout} follows from {\sc {\cblue
identity}}+{\sc coproduct}. We obviously need to consider only
$\gl
>1$ and then we use the fact that given \mor s $h_i:A\ra B_i\,
(i<\gl)$, their wide pushout $h:A\ra C$ can be obtained from $
{\coprod_{i<\gl}h_i}$ by pushing out along the codiagonal $\nabla:
{\coprod_{\gl}A\ra A}$:

$$
\xy (0,0)*{\xymatrix{\coprod A\ar[rr]^{\coprod h_i}\ar[d]_{\nabla}&& \coprod B_i\ar[d]^{}\\
A \ar[rr]_{h}&& C}}="D"; (24,-14)*{}="A"; (24,-9)*{}="B";
(29.5,-9)*{}="C"; "A"; "B" **\dir{-}; "B"; "C" **\dir{-};
\endxy
$$

}\end{sub}

\vspace*{3mm}
\begin{sub}\label{aa2.19}{\em {\rem} The deduction system of the last lemma has five rules, but the advantage
 against the system \ref{aa2.11}
is that they are particularly simple to formulate:

\vspace*{4.5mm}

\hspace*{2cm}\begin{tabular}{p{2.7cm}l}{\cblue {\sc identity}}\\ \\
\end{tabular}\hspace*{3mm}\begin{tabular}{c}\hline $\id_A$\\ \end{tabular}

\hspace*{2cm}\begin{tabular}{p{2.7cm}l}{\sc
cancellation}\\ \\
\end{tabular}\hspace*{3mm}\begin{tabular}{c}$h_2  \cdot h_1$
\\ \hline $ h_1$ \\ \\ \end{tabular}

\hspace*{2cm}\begin{tabular}{p{2.7cm}l}{\sc
composition}\\ \\
\end{tabular}\hspace*{3mm}\begin{tabular}{c}$h_2\; \; h_1$
\\ \hline $ h_2\cdot h_1$ \\ \\ \end{tabular}
\begin{tabular}{l}\\ {\hspace*{10mm}if $h_2\cdot h_1$ is defined}\\ \\ \end{tabular}

\hspace*{2cm}\begin{tabular}{p{2.7cm}l}{\sc pushout}\\
\end{tabular}\hspace*{3mm}\begin{tabular}{c}$ h$\\ \hline
 $ h'$ \\ \end{tabular}
\begin{tabular}{l}\\ {\hspace*{10mm}given }\\ \\ \end{tabular}
\begin{tabular}{l}$
\xy (0,0)*{\xymatrix{\ar[r]^h\ar[d]& \ar[d]\\
\ar[r]^{h'}& }}="D"; (7,-10.5)*{}="A"; (7,-7)*{}="B";
(10.5,-7)*{}="C"; "A"; "B" **\dir{-}; "B"; "C" **\dir{-};
\endxy
$\end{tabular} \vspace*{3mm}

\hspace*{2cm}\begin{tabular}{p{2.7cm}l}{\sc coproduct}
\\ \end{tabular}\hspace*{3mm}\begin{tabular}{c} $h_i\; (i\in I)$ \\ \hline
 $ {\coprod_{i\in I} h_i}$\\ \end{tabular}

\vspace*{6mm}

We prove below that \ref{aa2.11}, and therefore the above
equivalent deduction system, is not only sound but (in a number of
categories) also complete.}\end{sub}

\bs{\em {\rem} To relate our deduction rules to the usual ones (of
classical logic), let us consider, as in the Introduction, the \cat
of all $\Sigma$-structures. Then any object $A$ can be presented by
a set $A^\m(X)$ of atomic formulas with parameters $X$ in $A$: for
the familiar algebraic structures, this is just the usual concept of
generators and relations. Given a morphism $f:A\ra B$, and such
presentations $A^\m(X)$ and $B_o^\m(Y)$ of $A$ and $B$, we can also
present $B$ by $B^\m(X,Y)$, which is the union of $B_o^\m(Y)$ and
the set of all the equations $x = t(Y)$ for which $f(x) = t(Y)$ ($t$
a $\Sigma$-term). Then for the sentence
$$f^\m:=\forall X(\wedge A^\m(X)\ra \exists Y(\wedge B^\m(X,Y)))$$
we have that an object $C$ is $f$-injective iff $C\models f^\m$.
Note that if $f$ is finitary (see the Introduction or 3.4 below),
the presentations, and hence $f^\m$, can be chosen to be finitary
(more details in \cite{AR}, 5.33). Now, we can associate
Gentzen-style rules to sets of atomic formulas, generalizing the
idea of what was done (with more accuracy) in \cite{ASS} for sets of
equations: associating
$$A^\m(X)\Rightarrow B^\m(X,Y)$$
to $\forall X(\wedge A^\m(X)\ra \exists Y(\wedge B^\m(X,Y)))$, the
{\sc identity} axiom is of course \vspace*{-4mm}
$$\begin{array}{c}\\ \hline A^\m(X)\Rightarrow A^\m(X)\end{array}\; \; ;$$
{\sc cancellation} is a categorical version of the ``restriction"
rule
$$\begin{array}{c}A^\m(X)\Rightarrow (B^\m(X,Y) \cup C^\m(X,Y,Z))\\\hline A^\m(X)\Rightarrow B^\m(X,Y)\end{array}\; \; ;$$
{\sc pushout} is essentially the ``weakening" rule
$$\begin{array}{c} A^\m(X)\Rightarrow B^\m(X,Y)\\ \hline (A^\m(X)\cup C^\m(X,Z))\Rightarrow B^\m(X,Y)\end{array}\; \; ;$$
and {\sc composition} is a `` cut" rule
$$\begin{array}{c}A^\m(X)\Rightarrow B^\m(X,Y),\; B^\m(X,Y)\Rightarrow C^\m(X,Y,Z)\\\hline
A^\m(X)\Rightarrow C^\m(X,Y,Z)\end{array}\; \; .$$ The usual
stronger ``cut" rule
$$\begin{array}{c}A^\m(X)\Rightarrow B^\m(X,Y),\; ((B^\m(X,Y) \cup C^\m(X,Y,Z))\Rightarrow D^\m(X,Y,Z,U)
\\\hline (A^\m(X)\cup C^\m(X,Y,Z))\Rightarrow D^\m(X,Y,Z,U)\end{array}$$
corresponds to $$\begin{array}{c}\xymatrix{A\ar[r]^f&B},\;
\xymatrix{B+C\ar[r]^<<<<<g &D}
\\\hline \xymatrix{A+C\ar[rr]^>>>>>>>>>>{g\cdot (f+1_C)}&& C}\end{array}\; \; ,$$
which is proved via

{ \begin{center}\begin{tabular}{p{4cm}}
\begin{tabular}{c}$\quad \; f \qquad \qquad g$\end{tabular}

 \begin{tabular}{p{3.2cm}}\hline \\ \end{tabular}

\vspace*{-4mm}

  \begin{tabular}{c}$\; \;\;f+\id_{C} \qquad g$\end{tabular}

  \begin{tabular}{p{3.2cm}}\hline \\ \end{tabular}

\vspace*{-4mm}

\begin{tabular}{c}$\quad \;\, g\cdot(f+1_C)$
\end{tabular}
\end{tabular}
\begin{tabular}{p{3.2cm}} \vspace*{-5mm} {\sc pushout}\\
 {\sc composition}
\end{tabular}

\end{center}}}\es

\section{Completeness in locally presentable categories}

\bs\label{b3.1}{\em \textbf{Assumption} In the present section we
study injectivity  in a \textit{locally presentable} \cat $\A$ of
Gabriel and Ulmer, see \cite{GU}
 or \cite{AR}.
 This means  that:
 \begin{enumerate}
    \item[(a)] $\A$ is cocomplete,
 \end{enumerate}
 and
  \begin{enumerate}
    \item[(b)]there exists {\cbb a regular} cardinal $\gl$ such that
    $\A$ has a set of $\gl$-presentable objects whose closure
    under $\gl$-filtered colimits is all of $\A$.
 \end{enumerate}
 Recall that an object $A$ is \textit{$\gl$-presentable} if its
 hom-functor hom($A,-):\A\ra \Set$ preserves $\gl$-filtered colimits. That is, given a $\gl$-filtered diagram $D$ with
 a colimit $c_i:D_i\ra C$ $(i\in I)$ in $\A$, then for every
 {\mor}
 $f:A\ra C$
 \begin{enumerate}
    \item[(i)] a factorization of $f$ through   $c_i$ exists for some $i\in
    I$,
 \end{enumerate}
 and
  \begin{enumerate}
    \item[(ii)] factorizations are essentially unique, i.e., given
    $i\in I$ and
    $c_i\cdot g^\m=c_i\cdot g^{\m\m}$
     for some
    $g^\m, g^{\m\m}:A\ra D_i$,
     there exists a connecting {\mor}
    $d_{ij}:D_i\ra D_j$
     of the diagram with
     $d_{ij}\cdot g^\m=d_{ij} \cdot g^{\m\m}$.
 \end{enumerate}
 }

 \es

 \bs\label{b3.2}{\em {\exas} (see \cite{AR}) Sets, presheaves, varieties of
 algebras and simplicial sets are examples of locally presentable
 categories. Categories such as $\mathbf{Top}$ (topological spaces) or
 $\mathbf{Haus}$ (Hausdorff spaces) are not locally
 presentable.}\es

 \bs\label{b3.3}{\em {\rem}
 (a) In the present section we prove that the \Ij Logic is complete in every locally presentable category.
The
 reader may decide to skip this section since
 we prove a more general result in Section 6. Both of our proofs are based on
 the fact
 that for every set $\ch$ of \mor s  the full sub\cat Inj$\ch$ (of all objects \ijw \mor s of $\ch$) is
  \textit{weakly reflective}.
 That is: every object $A\in \A$ has a morphism $r:A\ra \overline{A}$, called a weak reflection, such that
\begin{enumerate}
    \item[(i)] $\overline{A}$ lies in Inj$\ch$
\end{enumerate}
and
\begin{enumerate} \item[(ii)] every morphism from $A$ to an object of Inj$\ch$
 factors through $r$ (not necessarily uniquely).
\end{enumerate}
In the present section we will utilize the classical \textit{Small
Object Argument} of D. Quillen \cite{Q}: this tells us that every
object $A$ has a weak reflection $r:A\ra \overline{A}$ in Inj$\ch$
such that $r$ is a transfinite composite of \mor s of the class
$$\wth=\{ k; \, k \mbox{ is a pushout of a member of $\ch$ along some
morphism}\}.$$

(b) The reason for proving the completeness based on the Small
Object Argument in the present section is that the proof is short
and elegant. However, by using a more refined construction of weak
reflection in Inj$\ch$, which we present in Section 5, we will be
able to prove the completeness in the so-called strongly locally
ranked categories, which include $\mathbf{Top}$ and $\mathbf{Haus}$.

The spirits of the two proofs are  quite different. Given an
injectivity consequence $h$ of a set of morphisms, in this section
we will show how to derive a formal proof of $h$ from Quillen's
construction of the weak reflection; this construction is ``linear",
forming a transfinite composite. In the next section, a weak
reflection will be constructed as a colimit of a filtered diagram
which somehow presents simultaneously all the possible formal
proofs.}\es

\bs\label{b3.4}{\em {\defn} A morphism is called
\textit{$\gl$-ary} provided that its domain and codomain are
$\gl$-presentable objects. For $\gl=\aleph_0$ we say
\textit{finitary}. }\es

\bs\label{b3.4}{\em {\rem} (a) The $\gl$-ary morphisms are
precisely the $\gl$-presentable objects of the arrow category
$\A^\rightarrow$. In contrast, M. H\'ebert introduced in \cite{H}
$\gl$-presentable \mor s; these are the \mor s $f:A\ra B$ which
are $\gl$-presentable objects of the slice \cat $A\downarrow \A$.
In the present paper we will not use the latter concept.

(b) We work now with the \textit{Finitary \Ij Logic}, i.e., the
deduction system \ref{2.4} applied to finitary \mor s. We
generalize this to the $k$-ary logic below. }\es

 \bs\label{b3.6}{\em {\theo}} The Finitary \Ij Logic is
complete in every locally presentable category $\ca$. That is, given
a set $\ch$ of finitary \mor s in $\ca$, then every finitary
morphism $h$ which is an injectivity consequence of $\ch$ is a
formal consequence in the deduction system \ref{2.4}. Briefly:
$$\ch \models h \mbox{ implies } \ch\vdash h.$$ \es

 {\pf}
  Given a  finitary morphism $h:A\ra B$
which  is an injectivity consequence of $\ch$,
 we prove that $$\ch \vdash h.$$

 (a) {\cbb The} above
 object $A$ has a weak reflection
 $$r:A\ra \overline{A}$$ such that $r$ is a transfinite
 composition of \mor s in
 $\wth$,
  see \ref{b3.3}(a). Since $\ch\models h$, it  follows that
 $\overline{A}$ is \ijw $h$, which yields a {\mor} $u$
 forming a commutative triangle
 $$\xymatrix{A\ar[rr]^r\ar[dr]_h&&\ova\\&B\ar[ru]_u&}$$

 (b) Consider all commutative triangles as above where
 $r:A\ra \overline{A}$ is any
 $\ga$-composite of \mor s in $\wth$ for some ordinal $\ga$ and $u$ is arbitrary. We
 prove that the least possible $\ga$ is finite. This finishes the
 proof of $\ch \vdash h$: In case $\ga=0$, we have that $\id=u\cdot
 h$,
and we derive $h$ via {\sc identity} and {\sc cancellation}. In case
$\ga$ is a finite ordinal greater than 0, we have that $r$ is
provable from $\mathcal{H}$ using {\sc pushout} and {\sc
composition}. Consequently, via {\sc cancellation}, we get $h$.

Let $\mathcal{C}$ be the class of all ordinals $\ga$ such that there
are an $\ga$-composite $r$ of \mor s of $\wth$ and a morphism $u$
with $r=u\cdot h$. To show that the least member $\gamma$ of
$\mathcal{C}$ is finite, we prove that for each ordinal $\gamma\geq
\omega$ in $\mathcal{C}$ we can find another ordinal in
$\mathcal{C}$ which is smaller than $\gamma$.

 \vspace*{2mm}

\h A. Case $\gamma=\beta+m$, with $\beta$ a limit ordinal and $m>0$
finite.  Let $a_{i,i+1}\; (i< \beta+m)$ be  the corresponding chain
with $r=a_{0, \gb+m}$.
   Since $a_{\gb,
\gb+1}$ lies in $\wth$, we can express it as a pushout of some
{\mor} $k:D\ra D^\m$ in $\ch$:

 $\hspace*{3.5cm}\xy
(-70,0)*{\xymatrix{&&&&&&&D\ar[lllld]_q\ar[r]^k\ar[d]_p&D^\m\ar[d]^{p^\m}\\A_0\ar[r]_{a_{01}}&A_1\ar[r]_{a_{12}}
&&A_i\ar[r]_{a_{i,i+1}}
\ar[d]_{v_i}&A_{i+1}\ar[r]_>>>>{a_{i+1,i+2}}
\ar[d]_{v_{i+1}}&A_{i+2}\ar[r]_>>>>{a_{i+2,i+3}}
\ar[d]_{v_{i+2}}&&A_\gb\ar[r]_{a_{\gb,\gb+1}}
\ar[d]_{v_\gb}&A_{\gb+1}\ar@{=}[ld]\\&&&P_i\ar[r]_{p_{i,i+1}}&P_{i+1}\ar[r]_>>>>{p_{i+1,i+2}}
&P_{i+2}\ar[r]_>>>>{p_{i+2,i+3}}&&P_\gb&}}="D";
(-26,-14.2)*{}="A"; (-26,-9)*{}="B"; (-18.5,-9)*{}="C"; "A"; "B"
**\dir{-}; "B"; "C" **\dir{-};  ;
(-43.5,-14)*{\dots}; (-43.5,-29)*{\dots};
(-106.5,-14)*{\dots}\endxy$

\vspace*{0.7cm}

We have a colimit $A_\gb=\mbox{co}\hspace*{-0.7mm}\lim_{i<\beta}A_i$
of a chain of \mor s.   Hence, because $D$ is finitely presentable,
$p$ factorizes as $p=a_{i\gb}\cdot q$ for some $i<\gb$ and some
morphism $q:D\ra A_i$. Let $v_i$ be a \ps of $k$ along $q$, and form
a sequence $v_j$ of pushouts of $k$ along $a_{ij}\cdot q\,
(j<\beta)$ as illustrated in the diagram above (taking  colimits at
the limit ordinals). Then it is easily seen, due to
$p=a_{i\beta}\cdot q$, that $v_\gb
 =
\mbox{co}\hspace*{-0.7mm}\lim_{j<\beta}v_j$ is a pushout of $k$
along $p$. Thus, without loss of generality,
$$P_\gb= A_{\beta+1} \; \; \mbox{and}\; \; v_\gb=a_{\gb,\gb+1}.$$
Observe that, since $a_{j, j+1}$ lies in $\wth$, {\sc pushout}
implies that
$$p_{j,j+1}\in \wth \; \; \mbox{for all} \, \,  i\leq j
< \gb.$$ Also $v_i\in \wth$ since it is a \ps of $k$ along $q$.
Consequently, $a_{0,\beta+1}$ is a $\gb$-composite of  morphisms
$b_{j,j+1}\; (j< \beta)$   of $\wth$ as follows (where $l$ is the
first limit ordinal after $i$):
$$\begin{array}{lll}&b_{j,j+1}=a_{j,j+1}&\; \mbox{ for all $j<i$},\\
&b_{i,i+1}=v_i,\\
& b_{j,j+1}=p_{j-1,j}&\; \mbox{ for all $i<j<l$},\\
\mbox{and}&&\\
& b_{j,j+1}=p_{j,j+1}&\; \mbox{ for all $l\leq j<\gb$}.
\end{array}$$

Thus $r=a_{0,\beta+m}$ is a $(\beta+(m-1))$-composite of morphisms
of $\wth$.

\vspace*{2mm}

\h B. Case $\gamma$ is a limit ordinal.  The {\mor}
$$u:B\ra \ova=\mbox{co}\hspace*{-0.7mm}\lim_{i<\gamma}A_i$$ factors,
since $B$ {\cbb is finitely presentable}, through some $a_{i\gamma},
\; i<\gamma$:
$$u=a_{i\gamma} \cdot \overline{u}\; \; \mbox{for some $\overline{u}:B\ra A_i$.}$$ The parallel pair
$$\xymatrix{A=A_0 \ar@<.8ex>[r]^{\; \; \; \; \; \; \overline{u}\cdot h}\ar@<-.8ex>[r]_{\; \; \; \; \; \; a_{0i}}&A_i}$$
 is clearly merged
by the colimit {\mor} $a_{i\gamma}$ of
$A_\gamma=\mbox{co}\hspace*{-0.7mm}\lim_{i<\gamma}A_i$. Since $A$ is
 finitely presentable, hom($A,-$) preserves that colimit,
consequently (see (ii) in \ref{b3.1}.b), the parallel pair is also
merged by a connecting {\mor} $a_{ij}:A_i\ra A_j$ for some
$i<j<\gamma$:
$$a_{ij}\cdot \overline{u}\cdot h=a_{0j}.$$
This gives us a commutative triangle

$$\xy (0,0)*{   \xymatrix{A_0\ar[r]^{a_{01}}\ar[dr]_h&A_1\ar[r]^{a_{12}}& & A_j\\
&B\ar[rru]_{a_{ij}\cdot \overline{u}}&&   }}="A";
(32.5,0)*{{\dots}} ="B";\endxy
$$

\vspace*{1.3cm} \hspace*{-\parindent}thus $a_{0j}$ is a
$j$-composite of \mor s of $\widehat{\ch}$ with $j<\gamma$.

\bs\label{b3.7}{\em {\rem}  The above theorem immediatly generalizes
to the $k$-ary \Ij Logic, i.e., to the deduction system of \ref{ver}
applied to $k$-ary \mor s. Recall that for every set of objects in a
locally presentable \cat there exists a cardinal $k$ such that all
these objects are $k$-presentable. Consequently, for every set
$\ch\cup \{h\}$ of \mor s there exists $k$ such that all members are
$k$-ary. The proof that $\ch \models h$ implies $\ch\vdash h$ is
completely analogously to \ref{b3.6}: We show that the least
possible $\ga$ is smaller than $k$, thus  in Cases A. and B. we work
with $\gamma\geq k$. } \es

\bs\label{b3.8}{\em {\cor}} The \Ij Logic is sound and complete in
every locally presentable category.
 \es

In fact, given $$\ch \models h$$ find a cardinal $k$ such that all
members of $\ch\cup \{h\}$ are $k$-ary \mor s. Then $h$ is a
formal consequence of $\ch$ by \ref{b3.7}.

\bs\label{b3.9}{\em {\rem} The above corollary also follows from the
Small Object Argument (see \ref{b3.3}(a)): if $h:A\ra B$ is an
injectivity consequence of $\ch$ and if $r:A\ra \ova$ is the
corresponding weak reflection, then $r$ is clearly a formal
consequence of $\ch$. Since $\ova$ is \ijw $h$, it follows that $r$
factors through $h$, thus, $h$ is a formal consequence of $r$ (via
{\sc cancellation}).}
 \es

\section{Strongly locally
ranked categories}

\bs \label{a3.1}{\em {\rem} Recall that a \textit{factorization
system} in a category is a pair $(\ce,\, \cm)$ of classes of \mor
s containing all iso\mor s and closed under composition such that
\begin{enumerate}
\item[(a)] every {\mor} $f:A\ra B$ has a factorization $f=m\cdot
e$ with $e:A\ra C$ in $\ce$ and $m:C\ra B$ in $\cm$
\end{enumerate}
and
\begin{enumerate}
\item[(b)] given another such factorization $f=m^\m\cdot e^\m$
there exists a unique ``diagonal fill-in" {\mor} $d$ making the
diagram
$$\xymatrix{A\ar[r]^{e}\ar[d]_{e^\m}&C\ar[ld]_{d}\ar[d]^{m}\\C^\m\ar[r]_{m^\m}&B}$$
commutative.

\end{enumerate}

The factorization system is called \textit{left-proper} if every
morphism of $\ce$ is an epi{\mor}. In that case the
$\ce$-quotients of an object $A$ are the quotient objects of $A$
represented by \mor s of $\ce$ with domain $A$.

}\es

\bs \label{a3.2}{\em {\defn} Let $(\ce,\cm)$ be a factorization
system. We say that an object $A$ has \textit{$\cm$-rank $\gl$},
where $\gl$ is a regular cardinal, provided, that

\begin{enumerate}
\item[(a)] hom($A,-$) preserves $\gl$-filtered colimits of diagrams of
$\cm$-morphisms (i.e., given a $\gl$-filtered diagram $D$ whose
connecting \mor s lie in $\cm$, then every {\mor} $f:A\ra
\mbox{colim}D$ factors, essentially uniquely, through a colimit
map of $D$)
\end{enumerate}
and
\begin{enumerate}
\item[(b)] $A$ has less than $\gl$ $\ce$-quotients.

\end{enumerate}

 If
$\gl=\aleph_0$ we say that the object $A$ has \textit{finite
$\cm$-rank}.} \es

\bs \label{a3.3} {\em {\exas}   (1) For the factorization system
(Iso, All), rank $\gl$ is equivalent to  $\gl$-presentability.

(2) In the \cat $\mathbf{Top}$ of topological spaces, choose
$(\ce, \, \cm)$ = (Epi, Strong Mono). Here the $\cm$-subobjects
are  precisely the embeddings of subspaces. Every topological
space $A$ of cardinality $\alpha$ has $\cm$-rank $\gl$ whenever
$\gl>2^{2^\alpha}$. In fact,  hom($A,-$) preserves $\gl$-directed
unions of subspaces since $\alpha<\gl$. And the amount of quotient
objects of $A$ (carried by epimorphisms) is at most
$\sum_{\beta\leq \alpha}E_\beta T_\beta$ where $E_\beta$ is the
number of equivalence relations on $A$ of order $\beta$ and
$T_\beta$ is the number of topologies on a  set of cardinality
$\beta$. Since $E_\beta$ and $T_\beta$ are both  $\leq
2^{2^\beta}$,  we have {\cbb $\sum_{\beta\leq \alpha}E_\beta
T_\beta \leq  \ga \cdot 2^{2^\ga}\cdot 2^{2^\ga}<\gl$, thus we
conclude that $A$ has less than $\gl$ quotients.}} \es

\bs \label{new4.4}{\em {\rem} Every $\ce$-quotient of an object of
$\cm$-rank $\gl$ also has $\cm$-rank $\gl$. In fact (a) in
\ref{a3.2} follows easily by diagonal fill-in, and (b) is
obvious.}\es

\bs \label{a3.4} {\em {\defn}  A \cat $\A$ is called
\textit{strongly locally ranked} provided that it has a left-proper
factorization system $(\ce, \, \cm)$ such that
\begin{enumerate}
    \item[(i)] $\A$ is cocomplete;
    \item[(ii)] every object has an $\cm$-rank, and all objects of
    the same  $\cm$-rank form a set up to isomorphism;
    \item[(iii)] for every cardinal $\mu$ the collection of all objects of $\cm$-rank $\mu$
    is closed under $\ce$-quotients and
    under $\mu$-small colimits,
    i.e., colimits of diagrams with
    less than $\mu$ morphisms;
   \end{enumerate}  \h and
   \begin{enumerate} \item[(iv)] the sub\cat of all objects of $\A$ and all \mor s
    of $\cm$ is closed under filtered colimits in $\A$. \end{enumerate}

{\rem} The statement (iv) means that, given a filtered colimit
with connecting morphisms in $\cm$, then
\begin{enumerate}
    \item[(a)] the colimit cocone  is formed by \mor s of $\cm$
    \end{enumerate}
    \h and
    \begin{enumerate}
    \item[(b)] every other cocone of $\cm$-\mor s has the unique
    factorizing morphism in $\cm$.

\end{enumerate}

} \es

\bs \label{a3.5} {\em {\exas} (1) Every locally presentable category
is strongly locally ranked: choose
$$\ce\equiv \mbox{iso\mor s}, \, \cm\equiv \mbox{all \mor s.}$$
 In fact, see \cite{AR}, 1.9 for the
proof of (ii), whereas (iii) and (iv) hold trivially.

(2) Choose
$$\ce\equiv \mbox{epi\mor s}, \, \cm\equiv \mbox{strong mono\mor s.}$$
Here categories such as $\mathbf{Top}$ (which are not locally
presentable) are included. In fact, for a space $A$ of cardinality
$\alpha$ we have that hom($A,-$) preserves $\gl$-filtered colimits
(=unions) of subspaces whenever $\alpha<\gl$. Thus, by choosing a
cardinal $\gl>\alpha$ bigger than the number of quotients of $A$
we get an $\cm$-rank of $A$. It is easy to verify (iii) and (iv)
in $\mathbf{Top}$.

(3) Let $\cb$ be a full, isomorphism closed, $\ce$-reflective
sub\cat of a strongly locally ranked \cat $\A$. If  $\cb$ is closed
under filtered colimits of $\cm$-\mor s in $\A$, then $\cb$ is
strongly locally ranked. In fact, $\cb$ is closed under $\cm$ in the
sense that given $m:A\ra B$ in $\cm$ with $B\in \cb$, then $A\in
\cb$. (Indeed, we have a reflection $r_A:A\ra A^\m$ in $\ce$ and
$m=m^\m \cdot r_A$ for a unique $m^\m$; this implies that $r_A\in
\ce$ is an isomorphism, thus, $A\in \cb$.) Therefore the restriction
of $(\ce,\, \cm)$ to $\cb$ yields a factorization system. It fulfils
(ii)-(iv) of \ref{a3.4} because $\cb$ is closed under filtered
colimits of $\cm$-morphisms.

(4) The \cat $\mathbf{Haus}$ of Hausdorff spaces is strongly locally
ranked: it is an epireflective sub\cat of $\mathbf{Top}$ closed
under filtered unions of subspaces.

}\es

\bs \label{a3.6} {\em \textbf{Observation} In a strongly locally
ranked \cat the class $\cm$ is closed under transfinite composition.
This follows from (iv).} \es

\bs \label{a3.7}{\em {\defn} A {\mor} is called \textit{$k$-ary}
if its domain and codomain have $\cm$-rank $k$. In case $k
={\aleph_0}$ we speak of \textit{finitary \mor s.}} \es

\bs \label{a3.8}{\em {\rem} The name ``strongly locally ranked" was
chosen since our requirements are somewhat stronger than those of
\cite{AHRT}: there a \cat is called locally ranked in case it is
cocomplete, has an $(\ce,\cm)$-factorization, is $\ce$-cowellpowered
and for every object $A$ there exists an infinite cardinal $\gl$
such that hom$(A,-)$ preserves colimits of $\gl$-chains of
$\cm$-monomorphisms.  Our definition of rank and the condition
\ref{a3.4}(ii) imply that the given category is $\ce$-cowellpowered.
Thus, every strongly locally ranked category is locally ranked.

An example of a locally ranked category that is not strongly locally
ranked is the category of $\sigma$-semilattices (posets with
countable joins and functions preserving them): condition
\ref{a3.4}(iv) fails here. Consider e.g. the $\omega$-chain of the
posets exp$(n)$ (where $n=\{0,1,\dots,n-1\}$), $n\in \omega$,  with
inclusion as order. The colimit of this chain is exp$(\mathbb{N})$
ordered by inclusion. If $M$ is the poset of all finite subsets of
$\mathbb{N}$ with an added top element, then the embeddings
exp$(n)\hookrightarrow M$ form a cocone of the chain, but the
factorization morphism exp$(\mathbb{N}) \ra M$ is not a
monomorphism.} \es

\section{A construction of weak reflections}

\bs \label{a4.1} {\em {\textbf{Assumption}} In the present section
$\A$ denotes a strongly locally ranked category. For every infinite
cardinal $k$, $\A_k$ denotes a chosen  set of objects of $\cm$-rank
$k$ closed under $\ce$-quotients and $k$-small colimits. In
particular, one may of course choose $\A_k$ to be a set of
representatives of all the objects of $\cm$-rank $k$ up to
isomorphism.}\es

Given  a set $\ch \subseteq \cm$ of $k$-ary \mor s of $\A_k$
(considered as a full sub\cat of $\A$), \cite{AHRT} provides a
construction of a weak reflection   in Inj$\,\ch$, which
generalizes the Small Object Argument (see \ref{b3.3}). However,
this does not appear to be sufficient to prove our Completeness
Theorem for the finitary case. The aim of this section is to
present a different, more appropriate construction.

We begin with the case $k=\omega$ and come back to the general case
at the end of this section.

\bs \label{a4.2} {\em \textbf{Convention} (a) Morphisms with
domain and codomain in $\A_\omega$ are called \textit{petty}.

(b) Given a set $\ch$ of petty \mor s, let
$$\overline{\ch}$$denote the closure of $\ch$ under finite composition and pushout in $\A_\omega$. (That is,
$\ovh$ is the closure  of $\ch \cup \{\id_A;\, A \in \A_\omega\}$
under binary composition and pushout along petty \mor s.)

(c) Since $\overline{\ch} \subseteq \mbox{mor}\A_\omega$ is a set,
we can, for every object $B$ of $\A_\omega$, index all \mor s of
$\overline{\ch}$ with domain $B$ by a set -- and that indexing set
can be chosen to be independent of $B$. That is, we assume that a
set $T$ is given and that for every object $B\in \A_\omega$,
\begin{equation}\label{(t)}\{h_B(t):B\ra B(t)\; ; \; t\in T \}\end{equation}
is the set of all \mor s of $\overline{\ch}$ with domain $B$. }\es

\bs \label{a4.3} {\em \textbf{Diagram $\mathbf{D_A}$} For every
object $A\in \A_\omega$ we define  a diagram $D_A$ in $\A$ and
later prove that a weak reflection of $A$ in Inj$\,\ch$ is
obtained as a colimit of $D_A$. The domain $\cd$ of $D_A$,
independent of $A$, is the poset of all finite {\bl words}
$$\varepsilon,\, M_1,\, M_1M_2,\, \dots,\, M_1\dots M_k\; \; (k<\omega)$$
{\bl where $\varepsilon$ denotes the empty word and} each $M_i$ is
a finite subset of $T$. The ordering is as follows:

$$M_1\dots M_k\leq N_1\dots N_l\; \; \mbox{ iff } \; \; k\leq l \; \; \mbox{ and }\;
\; M_1\subseteq N_1, \, \dots,\, M_k\subseteq N_k.$$ Observe that
$\varepsilon$ is the least element.

We denote the objects $D_A(M_1\dots M_k)$ of the diagram ${D_A}$
by
$$A_M \mbox{ where } M=M_1\dots M_k,$$ and if $M_1 \dots M_k\leq
N_1 \dots N_l \,= N$, we denote by
$$a_{M,N}:A_M\ra A_N$$ the corresponding
connecting morphism of $D_A$. We
define these objects and connecting morphisms by induction on the
length $k$ of the word $M=M_1 \dots M_k$ considered.

\textit{Case $k=0$}: $A_{\varepsilon}=A$.

\textit{Induction step}: Assume that all objects $A_M$ with $M$ of
length less than or equal to $k$ and all connecting morphisms
between them are defined. For every word $M$ of length $k+1$
denote by
$$M^\star\leq M$$ the prefix of $M$ of length $k$, and  define
the object $A_M$ as a colimit of the following finite diagram
$$
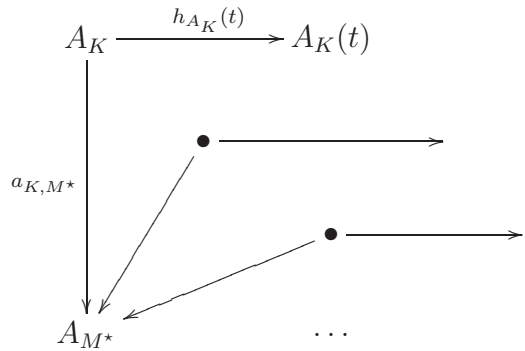
\xymatrix{A_K\ar[rr]^{h_{A_K}(t)}\ar[ddd]_{a_{K,M^\star}}&&A_K(t)&&\\&\bullet\ar[rr]\ar[ldd]&&&\\
&&\bullet\ar[rr]\ar[lld]&&\\A_{M^\star}&&\dots&&}$$ where $K$ ranges
over all words $K\in \mathcal{D}$ with $K\leq M^\star$ and $t$
ranges over the set $M_{k+1}$. Thus, $A_M$ is equipped with (the
universal cone of) \mor s
$$a_{M^\star,M}:A_{M^\star}\ra A_M \; \; \mbox{(connecting morphism of
$D_A$)}$$ and
$$d^K_M(t):A_K(t)\ra A_M \; \; \mbox{for all $K\leq M^\star,\, t\in M_{k+1}$,}$$
forming commutative squares
\begin{equation}\label{5.1a}\xymatrix{&
A_K\ar[ld]_{a_{K,M^\star}} \ar[rd]^{h_{A_K}(t)}&\\
A_{M^\star}\ar[rd]_{a_{M^\star,M}}&&A_K(t)\ar[ld]^{d_M^K(t)}\\&A_M&}\end{equation}
This defines the objects $A_M$ for all words of length $k+1$. Next
we define connecting \mor s
$$a_{N,M}:A_N\ra A_M$$
for all words $N\leq M$. If the length of $N$ is at most $k$, then
$N\leq M^\star$ and we define $a_{N,M}$ through the (already
defined) connecting morphism $a_{N,M^\star}$ by composing it with
the above $a_{M^\star,M}$. If $N$ has length $k+1$, we define
$a_{N,M}$ as the unique morphism for which the diagrams
\begin{equation}\label{5.1b}\xymatrix{&&&& A_K\ar[ld]_{a_{K,N^\star}} \ar[rd]^{h_{A_K}(t)}&&\\&&& A_{N^\star}\ar[rd]^{a_{N^\star,N}}
\ar[dddr]_{a_{N^\star, M}}&&A_K(t)\ar[ld]_{d_N^K(t)}
\ar[lddd]^{d^K_M(t)}&\\&&&&A_N\ar[dd]^<<<<<{a_{N,M}}&&(K\leq
N^\star,\, t\in N_{k+1})\\&&&&&\\&&&&A_M&&}\end{equation} commute.

It is easy to verify that the \mor s $a_{N,M}$ are well-defined
and that $D_A:\mathcal{D}\ra \A$ preserves composition and
identity \mor s.

}\end{sub}

\begin{sub}\label{3.3}\label{a4.6}{\em {\lem} }All connecting morphisms of the diagram $D_A$
lie in $\overline{\mathcal{H}}$.

{\em {\pf} We first observe that, given a finite diagram
$$\xymatrix{A_i\ar[r]^{h_i}\ar[d]_{f_i}&B_i\\C&}\; \;  \,\; \begin{array}{l}\\ \\ \\(i\in I) \,\end{array}$$ with all
 $h_i$ in $\ovh$,   a colimit
\begin{equation}\label{d1}\xymatrix{A_i\ar[r]^{h_i}\ar[d]_{f_i}&B_i\ar[d]^{d_i}\\C\ar[r]_{h}&D}\; \;  \,\; \begin{array}
{l}\\ \\ \\(i\in I) \,\end{array}\end{equation} is obtained  by
first considering  pushouts $h_i^\m$ of $h_i$ along $f_i$ and then
forming a wide pushout $h$ of all $h_i^\m\, (i\in I)$. Consequently,
the  connecting morphisms of $D_A$ are formed by repeating one of
the following steps: a finite wide pushout of \mor s in $\ovh$, a
composition of \mor s in $\ovh$, and a {\ps} of a {\mor} in $\ovh$
along a petty {\mor}. Since $\overline{\mathcal{H}}$ is closed, by
\ref{a4.2}, under the latter,  it is closed under the first one in
the obvious sense, see the construction of a finite wide pushout
described in Example \ref{2.7}.

}\es

\begin{sub}\label{3.4}\label{a4.7}{\em {\lem} }For every object $A_M$ of the diagram $D_A$ and every morphism $h:A_M\ra B$ of
 $\overline{\mathcal{H}}$ there exists a connecting morphism $a_{M,\, N}:A_M\ra A_N$ of $D_A$ which factors through $h$.

{\em {\pf} We have $M=M_1\dots M_k$ and $h=h_{A_M}(t)$ for some
$t\in T$. Put $$N=M_1\dots M_k\{t\}.$$ Then for $K=M$ the
definition of  $d_{N}^K(t)$ (see \eqref{5.1a}) gives  the
following commutative diagram:
$$\xymatrix{A_{M}\ar[rr]^{h_{A_{M}}(t)}\ar[d]_{id}&&A_{M}(t)\ar[d]^{d_N^K(t)}\\
A_{M}\ar[rr]_{a_{M,N}}&&A_{N}}$$ Consequently,
$$a_{M,N}=d_N^K(t)\cdot h_{A_M}(t)$$ as required. {\hfill
 }}\end{sub}

\begin{sub}\label{3.5}\label{a4.8}{\em {\prop} } Let $\ch$ be a set
of petty \mor s with $\ovh \subseteq \cm$. Then for every object $A \in \A_\omega$
a colimit $\gamma_M:A_M\ra \hat{A}\; (M\in\mathcal{D})$ of the diagram $D_A$
yields a weak reflection
 of $A$ in Inj$\, \mathcal{H}$ via
 $$r_A=\gamma_{\varepsilon}:A\ra \hat{A}.$$

{\em {\pf} (1) $\ha$ is \ijw $\ch$: We want to prove that given $h
\in\mathcal{H}$ and $f$ as follows
$$\xymatrix{B\ar[r]^h\ar[d]_f&C\\\hat{A}}$$
then $f$ factors through $h$. Firstly, since
$\hat{A}=\mbox{colim}D_A$ is a directed colimit of $\ovh$-\mor s
(see \ref{a4.6}) with $\ovh \subseteq \cm$, and $B$ has finite
$\cm$-rank (because $B\in \A_\omega$),
  it follows that hom($B,-$) preserves
the colimit of $D_A$. Thus,  there exists a colimit morphism
$\gamma_N:A_N\ra \hat{A}$ through which $f$ factors,
$f=\gamma_N\cdot f^\m$. $$\xymatrix{B\ar[r]^h\ar[d]_f\ar[dr]_{f^\m}&C\ar[dr]^{f^{\m\m}}&\\
\hat{A}&A_N\ar[l]^{\gamma_N}\ar[d]^{a_{N,\, M}}\ar[r]^{h^\m}&A_N(t)\ar[dl]^{h^{\m\m}}\\
& A_M\ar[ul]^{\gamma_M}&}$$ By pushing $h\in \mathcal{H}$ out
along $f^\m$ we obtain a morphism $h^\m \in
\overline{\mathcal{H}}$. Then by \ref{3.4} there exists $M\geq N$
such that $a_{N,M}=h^{\m\m}\cdot h^\m$ for some
$h^{\m\m}:A_N(t)\ra A_M$. The above commutative diagram proves
that $f$ factors through $h$.

(2) Let $B$ be \ijw$\mathcal{H}$. For every morphism $f:A\ra B$ we
define a compatible cocone $f_M:A_M\ra B$ of the diagram $D_A$ by
induction on
$$k=\mbox{the length of the word }M$$
such that $f_{\varepsilon}=f$. Then the desired factorization of
$f$ is obtained via the (unique) factorization $g:\hat{A}\ra B$
with $g\cdot \gamma_M=f_M$: in fact, $g\cdot r_A=f$.

For $k\mapsto k+1$, choose for every word $N$ of length $k$ and
every $t\in T$ a morphism $f_N(t)$ forming a commutative triangle
$$\xymatrix{A_N\ar[r]^{h_{A_N}(t)}\ar[d]_{f_N}&A_N(t)\ar[dl]^{f_N(t)}\\B&}$$
(recalling that $B$ is $\overline{\mathcal{H}}$-injective because
it is $\mathcal{H}$-injective). Then for every word $M$ of length
$k+1$ we  have a unique factorization $f_M:A_M\ra B$ making the
following diagrams
\begin{equation}\label{(3.5)}\xymatrix{A_{K}\ar[rr]^{h_{A_{K}}(t)}\ar[d]_{a_{K,\,
M^{\star}}}&&A_{K}(t)\ar[d]_{d_{M}^K(t)}
\ar[ddr]^{f_K(t)}&\\
A_{M^\star}\ar[drrr]_{f_{M^\star}}\ar[rr]^{a_{M^\star,M}}&&
A_M\ar[dr]_<<<{f_M}&\\&&&B}\end{equation} commutative for all
$K\leq M^\star$ and $t\in M_{k+1}$.

 {\bb Let us verify the compatibility
\begin{equation}\label{(3.7)}f_M=f_N\cdot a_{M,\, N}\qquad \mbox{ for all $M\leq N$ in $\mathcal{D}$.}\end{equation}
The last  diagram yields $f_{M^\star}=f_M\cdot a_{M^\star,\, M}$.
Therefore, it is sufficient to prove \eqref{(3.7)} for words $M$
and $N$ of the same length $k+1$. In order to do that, we will
show  that
\begin{equation}\label{*}f_M\cdot  d^K_{M}(t)=f_N\cdot a_{M,N}\cdot d^K_{M}(t),\,
\mbox{ for all $K\leq M^\star$ and $t\in M_{k+1}$,}
\end{equation}
and
\begin{equation}\label{**}f_M\cdot a_{M^\star,M}=f_N\cdot a_{M,N}\cdot
a_{M^\star,M}.
\end{equation}
Concerning \eqref{*}, we have
$$\begin{array}{ll}f_M\cdot
d^K_{M}(t)&=f_K(t)\\ \\&=f_N\cdot d^K_N(t), \; \mbox{by replacing $M$ by $N$ in \eqref{(3.5)}}\\
\\&=f_N\cdot a_{M,N}\cdot d^K_M(t),\; \mbox{by \eqref{5.1b}.}\end{array}$$
As for \eqref{**}, we have
$$\begin{array}{ll}f_M\cdot
a_{M^\star,M}&=f_{M^\star}\\ \\&=f_{N^\star}\cdot a_{M^\star,N^\star}\\ \\
&=f_N\cdot a_{N^\star,N}\cdot a_{M^\star, N^\star}, \; \mbox{by
replacing $M$ by $N$ in \eqref{(3.5)}}\\ \\ &=f_N\cdot
a_{M,N}\cdot a_{M^\star, M}.
\end{array}$$

}}\es

\bs \label{a4.9} {\em \textbf{Convention} Generalizing the above
construction from $\omega$ to any infinite cardinal $k$, we call
the \mor s of $\A_k$ \textit{$k$-petty}. Let us now denote by
$$\ovh_k$$the closure of $\ch$ under
$k$-composition (\ref{a2.9}) and \ps in $A_k$. Following
\ref{ver}, $\ovh_k$ is closed  under $k$-wide pushout. We again
assume that a set $T$ is given such that, for every object $B\in
\A_k$ we have an indexing $h_B(t):B\ra B(t)$, $t\in T$ of all \mor
s of $\ovh_k$ with domain $B$.}\es

\bs \label{a4.10} {\em \textbf{Diagram $\mathbf{D_A}$} The poset
$\cd$ of \ref{a4.3} is generalized to a poset $\cd_k$: Let
$\mathcal{P}_kT$ be the poset of all subsets of $T$ of cardinality
$<k$. The elements of $\mathcal{D}_k$ are all functions
$$M:\gl \ra \mathcal{P}_k T $$
where $\gl<k$ is an ordinal, including the case $\varepsilon:0\ra
\mathcal{P}_k T$. The ordering is as follows: for $N: \gl^{\m}\ra
\mathcal{P}_k T$ put
$$M\leq N \; \; \mbox{ iff }\; \; \gl\leq \gl^{\m}  \; \; \mbox{ and }\; \;
M_i\subseteq N_i \, \mbox{ for all $i<\gl$}.$$ We define, for
every $A\in \A_k$, the diagram $D_A:\mathcal{D}_k \ra
\mathcal{A}$. The objects $D_A(M)=A_M$  and the connecting \mor s
$a_{M,N}:A_M\ra A_N$ ($M\leq N$) are defined by transfinite
induction on $\gl<k$. For $\gl =0$ we have $A_{\varepsilon}=A$.
 The isolated step is precisely as in \ref{a4.3}, where for
$M:\gl+1\ra \mathcal{P}_k T$ we denote by $M^\star:\gl\ra
\mathcal{P}_k T$ the domain-restriction. The limit steps are defined
via colimits of smooth chains, see \ref{a2.9}: if $\gl<k$ is a limit
ordinal and $M:\gl\ra \mathcal{P}_k T$ is given, then $A_M$ is a
colimit of the chain $A_{M/i}\; (i<\gl)$, where $M/i$ is the domain
restriction of $M$ to $i$, with the connecting morphisms $a_{M/i,\,
M/j}:A_{M/i}\ra A_{M/j}$ for all $i\leq j<\gl$. The proof that these
chains are smooth is an easy transfinite induction.

It is also easy to see that all the above  results hold:
$\hat{A}=\mbox{colim} D_A$ is an $\mathcal{H}$-injective weak
reflection of $A$, and all connecting morphisms of $D_A$ are
members of $\overline{\mathcal{H}}$. Consequently, the proof of
the following proposition is analogous to that of \ref{a4.8}:}\es

\bs \label{5.9} {\em {\prop}} Let $\ch$ be a set of $k$-petty \mor
s with $\overline{\ch_k}\subseteq \cm$. Then for every object
$A\in \A_k$ a colimit $\gamma_M:A_M\ra \hat{A}$ of $D_A$ yields a
weak reflection of $A$ in Inj$\ch$ via
$r_A=\gamma_\varepsilon:A\ra\hat{A}$.\es

\section{Completeness in strongly locally ranked categories}

\bs \label{a5.1} {\em \textbf{Assumption} Throughout this section
$\ca$ denotes a strongly locally ranked category. We first prove the
completeness of the finitary logic. Recall that the finitary \mor s
are those where the domain and codomain are of finite $\cm$-rank.
Let us remark that whenever the class $\cm$ is closed under pushout,
then the method of proof of Theorem \ref{b3.6} applies again.
However, this excludes examples such as $\mathbf{Haus}$ (where
strong monomorphisms are not closed under pushout). }\es

\bs \label{a5.2} {\em {\theo}} The Finitary Injectivity Logic is
complete in every strongly locally ranked category. That is,  given
a set $\ch$ of finitary \mor s, every finitary morphism $h$ which is
an injectivity consequence of $\ch$ is a formal consequence (in the
deduction system of \ref{2.4}). Shortly: $\ch \models h$ implies
$\ch \vdash h$. \es

\bs \label{6.2a} {\em {\rem} We do not need the full strength of
weak local presentation for this result. We are going to prove the
completeness under the following milder assumptions on $\A$:
\begin{enumerate}
\item[(i)] $\A$ is cocomplete and has a left-proper factorization
system $(\ce,\, \cm)$; \item[(ii)] $\A_\omega$ is a set of objects
of finite $\cm$-rank, closed under finite colimits and
$\ce$-quotients; \item[(iii)] $\cm$ is closed under filtered
colimits in $\A$ (see \ref{a3.4} (iv)). \end{enumerate} The
statement we prove is, then, concerned with petty morphisms (see
\ref{a4.2}). We show that for every set $\ch$ of petty \mor s we
have
$$\ch \models h \, \;  \mbox{implies} \, \; \ch \vdash h  \; \; \mbox{(for all $h$ petty).}$$
The choice of $\A_\omega$ as a set of representatives of all
objects of finite $\cm$-rank yields the statement of the
theorem.}\es

 \h {\textbf {Proof of \ref{a5.2} and \ref{6.2a}}}   Let then $\ch$ be a set of petty \mor s, and let $$\ovh$$
 denote the closure of $\ch$ as in \ref{a4.2}.

 (1) We first prove that the theorem holds whenever $\ovh \subseteq \cm$. Moreover, we will show that
 for every petty injectivity consequence $\ch \models h$ we have a
 formal proof of
 $h$ from assumptions in $\ch$ such that the use of {\sc \ps} is always restricted to pushing out along petty \mor s.

  To prove
  this,
 consider, for the given petty injectivity consequence $h:A\ra B$ of $\ch$, the weak reflection $r_A:A\ra \ha$ in Inj$\,\ch$ of
 \ref{a4.8}. The object $\hat{A}$ is \ijw $h$, thus
$r_A$ factors through $h$ via some $f:B\ra \hat{A}$:
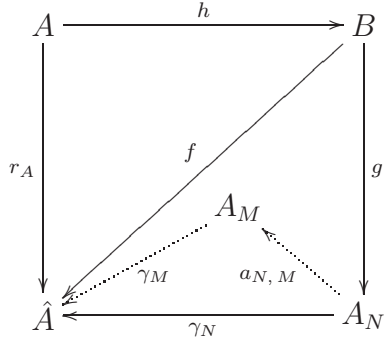
$$\xymatrix{A\ar[rrr]^h\ar[ddd]_{r_A}&&&B\ar[ddd]^g\ar[dddlll]_{f}\\
\\
&&A_M\ar@{.>}[dll]^{\gamma_M}&\\
\hat{A}&&&A_N\ar[lll]^{\gamma_N}\ar@{.>}[ul]^{a_{N,\, M}}}$$ Since
$B\in \A_\omega$, it has finite $\cm$-rank, and \ref{a4.6} implies
that hom($B,-$) preserves the colimit $\hat{A}=\mbox{colim} D_A$.
Then  $f$ factors through one of the colimit \mor s $\gamma_N:A_N\ra
\hat{A}$: $$f=\gamma_N\cdot g\, \; \mbox{ for some $g:B\ra A_N$}.$$
We know that $r_A=\gamma_{\varepsilon}$ is the composite of the
connecting morphism $a_{\varepsilon,\, N}:A\ra A_N$ of $D_A$ and
$\gamma_N$, therefore, $$\gamma_N\cdot a_{\varepsilon, \,
N}=r_A=\gamma_N\cdot g\cdot h.$$ That is, the colimit morphism
$\gamma_N$ merges the parallel pair $a_{\varepsilon,\, N},\, g\cdot
h:A\ra A_N$. Now the domain $A$ has finite $\cm$-rank, thus
hom($A,-$) also preserves $\ha=\mbox{colim}D_A$. Consequently, by
(ii) in \ref{b3.1}(b) the parallel pair is also merged by some
connecting morphism $a_{N,\, M}:A_N\ra A_M$ of $D_A$:
$$a_{N,\, M}\cdot a_{\varepsilon,\, N}=a_{N,\, M}\cdot g\cdot h:A\ra A_M.$$
The left-hand side is simply $a_{\varepsilon,\, M}$, and this is a
morphism of $\overline{\mathcal{H}}$, see Lemma \ref{a4.6}. Recall
that the definition of  $\overline{\mathcal{H}}$ implies that
every morphism in $\overline{\mathcal{H}}$ can be proved from
$\mathcal{H}$ using Finitary Injectivity Logic in which {\sc \ps}
is only applied to pushing out along petty \mor s. Thus, we have a
proof of the right-hand side $a_{N,\, M}\cdot g\cdot h$. The last
step is deriving $h$ from this by  {\sc cancellation}.

(2) Assuming $\ch \subseteq \ce$, then we prove that Inj$\,\ch$ is a
reflective sub\cat of $\A$, and for every object $A\in \A_\omega$
the reflection map $r_A:A\ra \hat{A}$ is a formal consequence of
$\ch$ lying in $\ce$:
$$\ch \vdash r_A \; \; \mbox{and}\; \; r_A \in \ce.$$ In fact,
from $\ch\subseteq \ce$ it follows that $\overline{\ch}\subseteq
{\ce}$ (since $\ce$ is closed under  composition and pushout).
Since $A$ has only finitely many $\ce$-quotients, see \ref{a3.2},
we can form a finite wide pushout, $r_A:A\ra \hat{A}$, of all
$\ce$-quotients of $A$ lying in $\overline{\ch}$. Clearly, $\ch
\vdash r_A$, in fact, $r_A \in \overline{\ch}$.

\h The object $\hat{A}$ is \ijw $\ch$: given $h:P\ra P^\m$ in
$\ch$ and $f:P\ra \hat{A}$, form a pushout $h^\m$ of $h$ along
$f$. This is an $\ce$-quotient in $\overline{\ch}$, then the same
is true for $h^\m\cdot r_A$. Consequently, $r_A$ factors through
$h^\m \cdot r_A$, and the factorization, $i:B\ra \hat{A}$, is an
epimorphism split by $h^\m$, thus, $f=i\cdot g\cdot h$:

$$\xymatrix{&P\ar[r]^{h}\ar[d]_{f}&P^\m\ar[d]^g\\
A\ar[r]^{r_A}&\ha \ar@<.8ex>[r]^{h^\m}&B\ar@<.8ex>[l]^{i}}$$

\h {\cbb The morphism} $r_A$ is a weak reflection: given a {\mor}
$u$ from $A$ to an object $C$ of Inj$\,\ch$, then $u$ factors
through $r_A$ because $C$ is \ijw $\ovh$ and $r_A\in
\overline{\ch}$.

(3) Let $\ch$ be arbitrary. We begin our proof by defining an
increasing sequence of sets  $\ce_i \subseteq \ce$ of petty \mor s
($i\in Ord$). For every member $f:A\ra B$ of $\overline{\ch}$ we
denote by $f_i$ a reflection of $f$ in Inj$\, \ce_i$:
$$\xymatrix{A\ar[r]^f\ar[d]_{r_A}&B\ar[d]^{r_B}\\
\ha\ar[r]_{f_i}&\hat{B}}$$

\textit{First step}: $\ce_0=\{\id_A; \, A\in \A_\omega\}$. Here
Inj$\, \ce_0=\A$, thus $f_0=f$.

\textit{Isolated step}:  For each $f\in \overline{\ch}$, let
$f_i=f_i^{\m\m}\cdot f_i^\m$ be the  $(\ce,\, \cm)$-factorization
of  the reflection $f_i$ of $f$ in Inj$\, \ce_i$, and put
$$\ce_{i+1}=\ce_i\cup\{ f_i^\m;\, f\in \overline{\ch}\}.$$

\textit{Limit step}: $\ce_j=\cup_{i<j} \ce_i$ for limit ordinals
$j$.

We prove that for every ordinal $i$ we have
\begin{equation}\label{fi'}\ch \vdash f^\m_i\; \; \; \mbox{for every }
f\in \overline{\ch}\end{equation} and
\begin{equation}\label{(*)}\mbox{Inj}\, \ch=\mbox{Inj}\, \ce_i\cap\mbox{Inj}\{f_i\}_{f\in
\overline{\ch}}.\end{equation}

{\cblue For $i=0$, \eqref{fi'} and \eqref{(*)} are trivial (use
{\sc cancellation} for \eqref{fi'} and {\sc identity} for
\eqref{(*)}). Given $i>0$, assuming that $\ch\vdash f^\m_j$ for
all $j<i$, with $f:A\ra B$ in $\ovh$, that is, $\ch \vdash \ce_i$,
we have, by (2), that
\begin{equation}\label{rB}\ch \vdash r_B\end{equation} where $r_B$ is the reflection
of $B$ in Inj$\, \ce_i$.} Thus, $\ch \vdash f_i\cdot r_A$.
Moreover, $r_A$ is an epi\mor , therefore the following square
$$\xymatrix{A\ar[r]^{r_A}\ar[d]_{r_A}&\ha\ar[r]^{f_i}&\hat{B}\ar[d]^{\id}\\
\ha\ar[rr]^{f_i}&&\hat{B}}$$ is a pushout, which  proves $\ch
\vdash f_i$ (via {\sc pushout}).  $\ch \vdash
 f^\m_i$ then follows by {\sc cancellation}.

 To prove \eqref{(*)}, observe that {\cblue \eqref{fi'} implies
 Inj$\,\ch\subseteq \mbox{Inj}\,  \ce_i$}, and our previous argument yields
  Inj$\,\ch \subseteq \mbox{Inj}\,\{f_i\}_{f\in \overline{\ch}}$. Thus, it remains to
 prove the reverse inclusion: every object $X$ \ijw
 $\ce_i\cup \{f_i\}_{f\in \overline{\ch}}
$ is \ijw $\ch$. In fact, given $f:A\ra B$ in $\ch$ and a {\mor}
$u:A\ra X$, then since $X\in \mbox{Inj}\, \ce_i$ we have a
factorization $u=v\cdot r_A$, and then the injectivity of $X$
w.r.t. $f_i$ yields the desired factorization of $u$ through $f$.
$$\xymatrix{A\ar[rr]^f \ar[dd]_{r_A}\ar[dr]^u&&B\ar[dd]^{r_B}\\ &X&&\\ \ha\ar[ur]^v\ar[rr]_{f_i}&& \hat{B}\ar@{-->}[ul]\\}$$

(4) Since $\A_\omega$ is a small category, there exists an ordinal
$j$ with $$\ce_j=\ce_{j+1}.$$ We want to apply (1) to the \cat
$$\A^\m=\mbox{Inj}\,\ce_j,$$
and the set $$\A^\m_\omega=\A_\omega\cap \mbox{obj}\A^\m.$$ Let us
verify that $\A^\m$ satisfies the assumptions (i) -- (iii) of
Remark \ref{6.2a} w.r.t. $$\ce^\m=\ce\cap \mbox{mor}\A^\m\; \;
\mbox{ and } \;  \; \cm^\m =\cm \cap \mbox{mor}\A^\m.$$

Ad(i): $\A^\m$ is cocomplete because it is reflective in $\A$.
Moreover, since the reflection maps lie in $\ce$, it follows that
$(\ce^\m,\, \cm^\m)$ is a factorization system: in fact, $\A^\m$
is closed under factorization in $\A$. Since $\ce \subseteq
\mbox{Epi}(\A)$, we have $\ce^\m \subseteq \mbox{Epi} (\A^\m)$.

Ad(iii): It is sufficient to prove that $\A^\m$ is closed under
filtered colimits of $\cm^\m$-\mor s in $\A$. In fact, let $D$ be
a filtered diagram in $\A^\m$ with connecting \mor s in $\cm$, and
let  $c_t:C_t\ra C \; (t\in T)$ be a colimit of $D$ in $\A$. Then
$C\in \A^\m$, i.e., $C$ is \ijw $f_j:\hat{A}\ra E$ for every $f\in
\overline{\ch}$. This follows from $\hat{A}$ having finite
$\cm$-rank (because $A\in \A_\omega$ implies $\hat{A}\in
\A_\omega$ due to the fact that $r_A:A\ra \hat{A}$ is an
$\ce$-quotient): since hom($\hat{A},-$) preserves the colimit of
$D$, every morphism $u:\hat{A}\ra C$ factors through some of the
colimit \mor s:
$$\xymatrix{&\ha\ar@{-->}[ld]_{v}\ar[r]^{f_j}\ar[d]^u&E\\C_t\ar[r]_{c_t}&C&
}$$ Since $C_t\in \A^\m$ is \ijw $f_j$, we have a factorization of
$v$ through $f_j$, and therefore, $u$ also factors through $f_j$.
This proves $C\in \A^\m$.

Ad(ii): Due to the above, every object of $\A^\m$ having a finite
$\cm$-rank in $\A$ has a finite $\cm^\m$-rank in $\A^\m$. Also, a
finite colimit of objects of $\A^\m$ in $\A^\m$ is a reflection
(thus, an $\ce$-quotient) of the corresponding finite colimit in
$\A$. Thus, it lies in $\A^\m_\omega$.

Next we claim that the set $\ch^\m=\{ f_j;\, f\in \ovh \}$ fulfils
$$\ch^\m  \subseteq \cm^\m$$
 and $\ch^\m$ is closed under petty identities, composition, and
pushouts along petty \mor s. In fact, in the above $(\ce,\,
\cm)$-factorization of $f_j$:
$$\xymatrix{A\ar[d]_{r_A}\ar[rr]^f&&B\ar[d]^{r_B}\\
\ha\ar[rr]_{f_j}\ar[rd]_{f^\m_j}&&\hat{B}\\&D\ar[ur]_{f^{\m\m}_j}
}$$ we know that $f_j^\m$ lies in $\ce_{j´+1}=\ce_j$ and $\ha$ is
\ijw $\ce_j$, thus, $f^\m_j$ is a split mono{\mor} (as well as an
epimorphism, since $\ce \subseteq \mbox{Epi}(\A)$). Thus, $f^\m_j$
is an iso\mor , which implies $f_j\in \cm$. $\ch^\m$ contains
$\id_A$ for every $A\in \A^\m_\omega$  because $\ovh$ contains it;
$\ch^\m$ is closed under composition because $\ovh$ is (and
$f\mapsto f_j$ is the action of the reflector functor from $\A$ to
Inj$\, \ce_j$). Finally, $\ch^\m$ is closed under pushout along
petty \mor s. In fact, to form a pushout of $f_j:\ha\ra \hat{B}$
along $u:\ha\ra C$ in $\A^\m=\mbox{Inj}\, \ce_j$, we form a
pushout, $g$, of $f$ along $u\cdot r_A$ in $\A$, and compose it
with the reflection map $r_D$ of the codomain $D$:
$$\xymatrix{&&A\ar[r]^f\ar[d]_{r_A}&B\ar[d]^{r_B}\ar[dddrr]^v&&\\
&&\ha\ar[r]^{f_j}\ar[ldld]_u&\hat{B}\ar[dr]_{\hat{v}}&&\\
&&&&\hat{D}&\\
C\ar[rrrrr]_g\ar[urrrr]^{\hat{g}}&&&&&D\ar[lu]^{r_D}}$$ Since $C$
lies in $\A^\m$, we can assume $r_C=\id_C$, and the reflection
$\hat{g}=r_D\cdot g$ of $g$ in $\A^\m$ is then a pushout of $f_j$
along $u$. Now $f\in \ovh$ implies $g\in \ovh$, and we have
$\hat{g}=g_j \in \ch^\m$.

(5) We are ready to prove that if a petty morphism $h:A\ra B$ is an
injectivity consequence of $\ch$, then $\ch\vdash h$ in $\A$. We
write $\ch\vdash_\A h$ for the latter since we work within two
categories: when we apply (1) to $\A^\m$  we use $\vdash_{\A^\m}$
for formal consequence in $\A^\m$. Analogously with $\models_{\A}$
and $\models_{\A^\m}$. Let $\hat{h}:\ha \ra \hat{B}$ be a reflection
of $h$ in $\A^\m$, then
$$\ch^\m\models_{\A^\m} \hat{h}$$
because every object $C\in \A^\m=\mbox{Inj}\,\ce_j$ which is \ijw
$\ch^\m=\{f_j\}_{f\in \ovh}$ is, due to \eqref{(*)}, \ijw $\ch$ in
$\A$. Then $C$ is \ijw $h$, and from $C\in \A^\m$ it follows
easily that $C$ is \ijw $\hat{h}$. Due to (4) we can apply (1).
Therefore,
$$\ch^\m\vdash_{\A^\m} \hat{h}.$$
We thus have a proof of $\hat{h}$ from $\ch^\m$ in $\A^\m$. We
modify it to obtain a proof of $h$ from $\ch$ in $\A$.  We have no
problems with a line of the given proof that uses one of the
assumptions $f_j\in \ch^\m$: we know from \eqref{fi'} that
$\ch\vdash_\A f_j$, and we substitute that line with a formal
proof of $f_j$ in $\A$. No problem is, of course, caused by the
lines using {\sc composition} or {\sc cancellation}. But we need
to modify the lines using {\sc pushout} because $\A^\m$ is not
closed under pushout in $\A$. However, a pushout, $g^{\m\m}$, of a
{\mor} $g$ along a petty {\mor} $u$ in $\A^\m$
$$\xymatrix{P\ar[r]^g\ar[d]_u&Q\ar[d]&\\
P^\m\ar[r]^{g^{\m\m}}\ar[drr]_{g^\m}&\\&&Q^\m\ar[ul]_{r_{Q^\m}}}$$
is obtained from a pushout, $g^\m$, of $g$ along $u$ in $\A$ by
composing it with a reflection map $r_{Q^\m}$ of the pushout
codomain. Recall that $P,\, P^\m,\, Q \in \A_\omega$ imply
$Q^\m\in \A_\omega$. Thus, we can replace the line $g^{\m\m}$ of
the given proof by using {\sc pushout} in $\A$ (deriving $g^\m$),
followed by a proof of $r_{Q^\m}$ (recall from  \eqref{rB} that
$\ch\vdash_\A r_{Q^\m}$) and an application of {\sc composition}.
We thus proved that $$\ch \vdash_\A \hat{h}.$$ Since $r_B\cdot
h=\hat{h}\cdot r_A$ and $\ch\vdash_\A r_A$  (see \eqref{rB}), we
conclude $\ch \vdash_\A \hat{h}\cdot r_A$; by {\sc cancellation}
then $\ch \vdash_\A h$.

\bs \label{compactness}{\em {\cor} (Compactness Theorem)} Let $\ch$
be a set of finitary \mor s in a strongly locally ranked category.
Every finitary morphism which is an injectivity consequence of $\ch$
is an injectivity consequence of a finite subset of $\ch$. \es

\bs {\em {\rem} We proceed by generalizing the completeness result
 from finitary to $k$-ary, where $k$ is an arbitrary infinite
cardinal. The \textit{k-ary logic}, then, deals with $k$-ary \mor
s (i.e., those having both domain and codomain of $\cm$-rank $k$)
and  the $k$-ary \Ij Deduction System of \ref{ver}. }\es

\bs \label{a6.2} {\em \theo} The $k$-ary \Ij Logic is complete in
every strongly locally ranked category. That is, given a set $\ch$
of $k$-ary \mor s, then every $k$-ary morphism which is an
injectivity consequence of $\ch$ is a formal consequence (in the
$k$-ary Injectivity Deduction System).

{\em {\pf} The whole proof is completely analogous to that of
Theorem \ref{a5.2}. As described in Remark \ref{6.2a} we work
under the following milder assumptions on the \cat $\A$:
\begin{enumerate}
    \item[(i)] $\A$ is cocomplete and has a left-proper factorization
    system $(\ce,\, \cm)$;
    \item[(ii)] $\A_k$ is a set of objects of $\cm$-rank $k$,
    closed under colimits of less than $k$ \mor s and under $\ce$-quotients;
    \item[(iii)] $\cm$ is closed under $k$-filtered colimits in
    $\A$.
\end{enumerate}
The statement we prove is concerned with $k$-petty \mor s (see
\ref{a4.9}). We denote by $\ovh_k$ the closure of {\cbb $\ch$ as
in \ref{a4.9}.} We write $\ch \vdash h$ for the $k$-ary \Ij Logic.

 (1) The theorem holds whenever
$\ovh_k\subseteq \cm$. The proof, based on the construction of a
weak reflection $\ha=\mbox{colim}D_A$ of \ref{a4.10}, is
completely analogous to that of (1) in \ref{a5.2}.

(2) Assuming $\ch\subseteq \ce$, then Inj$\,\ch$ is a reflective
subcategory, and the reflection maps $r_A$ fulfil  $\ch\vdash r_A$
and $r_A \in \ce$. This is analogous to the proof of (2) of
\ref{a5.2}.

(3) The definition of $\ce_i$ is precisely as in the proof of
\ref{a5.2}.

(4) For the first ordinal $j$ with $\ce_j=\ce_{j+1}$ the \cat
$\A^\m=\mbox{Inj}\, \ce_j$ fulfils the assumptions (i)-(iii)
above, and the set $\ch^\m=\{ f_j;\, f\in \ovh\}$ fulfils
$\ch^\m=\overline{\ch^\m} \subseteq \cm$.

(5) The theorem is then proved by applying (1) to $\A^\m$ and
$\ch^\m$: we get $\ch^\m \vdash \hat{h}$ in ${\A^\m}$ and we
derive $\ch\vdash h$ in $\A$ precisely as in the proof of
\ref{a5.2}.}\es

\bs \label{a6.3}{\em {\cor} } The \Ij Logic is sound and complete.
That is, given a set $\ch$ of \mor s of a strongly locally ranked
category, then the consequences  of $\ch$ are precisely the formal
consequences of $\ch$ (in the Injectivity Deduction System).
Shortly: $$\ch \models h \; \; \mbox{ iff } \; \; \ch \vdash h\; \;
\; \mbox{(for all \mor s $h$)}$$\es

In fact, soundness was proved in Section 2. Completeness follows
from Theorem \ref{a6.2}: since $\ch$ is a set, and since every
object of $\A$ has an $\cm$-rank, see \ref{a3.4}(ii), there exists
$k$ such that all domains and codomains of \mor s of $\ch \cup \{
h\}$ have $\cm$-rank $k$.

\section{Counterexamples}

\bs \label{a7.1} {\em {\exa} In ``nice" categories which are not
strongly locally ranked the completeness theorem can fail. Here we
refer to $\vdash$ of the Deduction System \ref{aa2.11} (and the
logic concerning arbitrary \mor s). We denote by
$$\mathbf{CPO(1)}$$ the \cat of unary algebras defined on $CPO$'s.
Recall that a $CPO$ is a poset with directed joins, and the
corresponding category, $\mathbf{CPO}$, has as \mor s the
\textit{continuous functions} (i.e., those preserving directed
joins). The \cat $\mathbf{CPO(1)}$ has as objects the triples
$(A,\sqsubseteq,\ga)$ where $(A, \sqsubseteq)$ is a $CPO$ and
$\ga:A\ra A$ is a unary operation. Morphisms are the continuous
algebra homo\mor s.

First let us observe that the assumption of cocompleteness is
fulfilled.} \es

\h {\lem} \textit{$\mathbf{CPO(1)}$ is cocomplete.}

{ {\pf} The \cat  $\mathbf{CPO}$ is   easily seen to be
cocomplete. The \cat $\mathbf{CPO(1)^*}$ of partial unary algebras
on $CPO$'s (defined as above except that we allow  $\ga:A^\m\ra A$
for any $A^\m\subseteq A$) is monotopological over $\mathbf{CPO}$,
see \cite{AHS}, since for every monosource
\newline{$f_i:(A,\sqsubseteq)\ra$} $(A_i,\sq_i,\ga_i)\; (i\in I)$ we
define a partial operation $\ga$ on $A$ at an element $x \in A$
iff $\ga_i$ is defined at $f_i(x)$ for every $i$, and then $$\ga
x=y\; \mbox{ iff } \; f_i(y)=\ga_i(f_i(x))\; \mbox{ for all $i\in
I$.}$$ Consequently, $\mathbf{CPO(1)^*}$ is cocomplete by
\cite{AHS}, 21.42 and 21.15. Further, $\mathbf{CPO(1)}$ is a full
reflective sub\cat of $\mathbf{CPO(1)^*}$: form a free unary
algebra on the given partial unary algebra, ignoring the ordering,
and then extend the ordering trivially (i.e., the new\ elements
are pairwise incomparable, and incomparable with any of the
original elements). Thus, $\mathbf{CPO(1)}$ is cocomplete.

\vspace*{2mm}

We will find \mor s $h_1,\, h_2$ and $k$ of $\mathbf{CPO(1)}$ with
$$\{h_1,\, h_2\}\models k\; \; \mbox{ but }\; \; \{h_1,\, h_2\}
\not\vdash k.$$

 (i) We define a  {\mor}
{$h_1$} that expresses, by injectivity, the condition

\h \begin{tabular}{lp{15cm}}(h1)\centerline{$x\sq\ga x \; \;
\mbox{ for all $x\in A.$}$}\end{tabular} \newline Let $=$ denote
the discrete order on the set $\mathbf{N}$ of natural numbers, and
$\sq$ that order enlarged by $0\sq 1$. Let $s:\mathbf{N}\ra
\mathbf{N}$ be the successor operation. Then
$$h_1=\id:(\mathbf{N},=,s)\ra (\mathbf{N},\sq, s)$$ is a {\mor} such that an algebra is \ijw $h_1$ iff it fulfils (h1) above.

(ii) The condition

\h
\begin{tabular}{lp{15cm}}(h2)\centerline{$A\not=\emptyset$}\end{tabular} is expressed by the injectivity w.r.t.
$$h_2:\emptyset\ra (\mathbf{N},=,s)$$ where $\emptyset$ is the empty (initial) algebra. The following {\mor} $k$
 expresses the existence
of a fixed point of $\ga$: $$k:\emptyset\ra 1$$ where 1 is a
one-element (terminal) algebra. }

\vspace*{2mm}

\h {{\prop}} {\it  $\{h_1,\, h_2\}\models k$ but $\{h_1,\,
h_2\}\not\vdash k$.}

{\pf} To prove $\{h_1,\, h_2\}\models k$, let $(A, \sq, \ga)$ be
\ijw $h_1$ and $h_2$, i.e., fulfill $x\sq \ga(x)$ and be nonempty.
Define a smooth (see \ref{a2.9}) chain $(a_i)_{i\in Ord}$ in
$(A,\sq)$ by transfinite induction: $a_0\in A$ is any chosen
element. Given $a_i$ put $a_{i+1}=\ga(a_i)$; we know that $a_i\sq
a_{i+1}$. Limit steps are given by (directed) joins,
$a_j=\bigsqcup_{i<j}a_i$. Since $A$ is small, there exist $i$ with
$a_i=a_{i+1}$, that is, $a_i$ is a fixed point of $\ga$. Thus, $A$
is \ijw $k$.

To prove  $\{h_1,\, h_2\}\not\vdash k$, it is sufficient to find
an extension $\mathcal{K}$ of the \cat $\mathbf{CPO(1)}$ in which
$\mathbf{CPO(1)}$ is closed under colimits (therefore $\vdash$ has
the same meaning in $\mathbf{CPO(1)}$ and in $\mathcal{K}$) and in
which there exists an object which is \ijw $h_1$ and $h_2$ but not
w.r.t. $k$. Thus $k$ cannot be proved in $\mathcal{K}$ from
$h_1,\, h_2$; consequently it cannot be proved in
$\mathbf{CPO(1)}$ either.

 We
define $\mathcal{K}$ by adding a single new object $K$ to
$\mathbf{CPO(1)}$. The only {\mor} with domain $K$ is $\id_K$. For
every algebra $(A,\sq,\ga)$ of $\mathbf{CPO(1)}$ we call a
function $f:A\ra Ord$ a \textit{coloring} of $A$ provided that it
is continuous and fulfils $f(\ga(x))=f(x)+1$ for all $x\in A$.

\h  The hom-object of $A$ and $K$ in $\mathcal{K}$ is defined to
be the class of all colorings of $A$. The composition in
$\mathcal{K}$ is defined ``naturally": given a continuous
homo{\mor}\newline $h:(A,\sq,\ga)\ra (B,\leq, \beta)$, then for
every coloring $f:B\ra Ord$ of $B$ we have a coloring $f\cdot
h:A\ra Ord$ of $A$. The \cat $\mathbf{CPO(1)}$ is a full sub\cat
of $\mathcal{K}$ closed under (small) colimits. In fact, given a
colimit cocone $a_i:A_i\ra A$ $(i\in I)$ in $\mathbf{CPO(1)}$,
then for every compatible cocone of colorings $f_i:A_i\ra Ord$
$(i\in I)$ there exists an ordinal $j$ such that all ordinals in
$\cup_{i\in I}f_i[A_i]$ are smaller than $j$. Let
$B=(j^+,\leq,\overline{s})$ be the object of $\mathbf{CPO(1)}$
where $\leq$ is the usual linear ordering of $j^+$ (the poset of
all ordinals smaller or equal to $j$), and $\overline{s}$ is the
successor map except $\overline{s}(j)=j$. Then the codomain
restriction $f^\m_i$ of each $f_i$ defines a continuous homo{\mor}
$f^\m_i:A_i\ra B$, and we obtain a compatible cocone
$(f_i^\m)_{i\in I}$ for our diagram. The unique continuous
homo{\mor} $g:A\ra B$ with $g\cdot a_i=f^\m_i$ yields, by
composing it with the inclusion $j^+\hookrightarrow  Ord$, a
coloring $f:A\ra Ord$ with $f\cdot a_i=f_i$ $(i\in I)$.

It is obvious that $K$ is \ijw $h_1$: every coloring of
$(\mathbf{N},=,s)$ is also a coloring of $(\mathbf{N},\sqsubseteq,
s)$. And $K$ is \ijw $h_2$ (because the inclusion
$\mathbf{N}\hookrightarrow Ord$ is a coloring of
$(\mathbf{N},=,s)$). But $K$ is not \ijw $k$, since $1$ has no
coloring.

\bs \label{a7.3} {\em {\exa}  None of the deduction rules of the
Finitary \Ij Deduction System   can be left out. For each of them
we present an example of a finite complete lattice $\A$ in which
the reduced deduction system is not complete (for finitary \mor
s).

 (1) {\cblue {\sc identity}}  The deduction
system {\sc cancellation}, {\sc composition} and  {\sc \ps} is not
complete because nothing can be derived from the empty set of
assumptions, although $\emptyset \models \id_A$.

(2) {\sc cancellation} In the poset
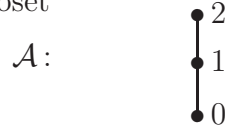
\begin{center}\begin{picture}(30,30)(0,0)
\put(-40,20){$\A\, $:} \put(30,0){\line(0,1){20}}
\put(27.1,-2){$\bullet\, 0$} \put(27.1,18){$\bullet\, 1$}
 \put(30,20){\line(0,1){20}}
 \put(27.1,36){$\bullet\, 2$}

\end{picture}
\end{center}
the only object \ijw $\{0\ra 2\}$ is $2$, thus, we see that
$\{0\ra 2\}\models 0\ra 1$.  However, $0\ra 1$ cannot be derived
from $0\ra 2$ by means of {\cblue {\sc identity}}, {\sc
composition} and {\sc \ps} because the set of all \mor s of $\A$
except $0\ra 1$ is closed under composition and pushout.

(3) {\sc composition} In $\A$ above we clearly have $\{ 0\ra 1,\,
1\ra 2\} \models 0\ra 2$. However, the set of all \mor s except
$0\ra 2$ is closed under left cancellation and pushout.

(4) {\sc \ps} In the poset
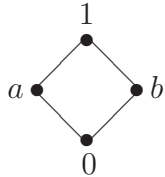
\begin{center}\begin{picture}(45,45)(0,-5)
\put(-40,20){$\; \; \; \;$} \put(12,20){\line(1,1){20}}
\put(29,-2){$\bullet$} \put(30,-12){0}
 \put(12,20){\line(1,-1){20}}
 \put(29,36){$\bullet$}
 \put(29,45){1}
 \put(52,20){\line(-1,1){20}}
 \put(52,20){\line(-1,-1){20}}
 \put(2,17){$a \, \bullet$}
 \put(48,17){$\bullet \, b$}

\end{picture}
\end{center}
we have $\{0\ra a\} \, \models \, b\ra 1$, but  we cannot derive
$b\ra 1$ from $0\ra a$ using {\cblue {\sc identity}}, {\sc
composition} and {\sc cancellation} because the set of all \mor s
except $b\ra 1$ is closed under composition and cancellation.}\es

\bs \label{a7.4} {\em {\exa} Here we demonstrate that in the
Finitary \Ij Logic we cannot restrict the statement of the
completeness theorem from the given strongly locally ranked category
$\A$ to its full subcategory $\A_\omega$ on all objects of finite
rank: although the relation $\vdash$ works entirely in $\A_\omega$,
the relation $\models$ does not.

More precisely, let $\ch \models_\omega h$ mean that every
$\ch$-injective object of finite $\cm$-rank is also $h$-injective.
And let $\vdash_\omega$ be the formal consequence w.r.t. Deduction
System \ref{2.4}. Then the implication

$$\ch \models_\omega h\;\mbox{ implies }\; \ch \vdash_\omega h$$
does NOT hold in general for sets  of finitary \mor s.

Indeed, let $\mathcal{A}=\mathcal{G}ra$ be the category of graphs,
i.e., binary relational structures $(A,R)$, $R\subseteq A\times
A$, and the usual graph homomorphisms. Recall that $\mathcal{G}ra$
is locally finitely presentable, and the finitely presentable
objects are precisely the finite graphs. Let us call a graph a
\textit{clique} if $R=A\times A-\Delta_A$. Denote by $C_n$ a
clique of cardinality $n$, and let $\mathbf{0}$ be the initial
object (empty graph).

For the set $$\mathcal{H}=\{\mathbf{0}\ra C_n\}_{n\in
\mathbb{N}}$$
 we have the following property:

 \centerline{every finite $\mathcal{H}$-injective graph $G$ has a loop (i.e., a morphism from $1$ to $G$).}
 \hspace*{-\parindent}In fact, if $G$ has cardinality
 less than $n$ and is injective w.r.t. $\mathbf{0}\ra C_n$, then we have a homomorphism $f:C_n \ra G$. Since $f$ cannot be one-to-one, there
 exist $x\not=y$ in $C_n$ with $f(x)=f(y)$ -- and the last element defines a loop of $G$ because $(x,y)$ is an edge of $C_n$.
 Hence
 $$\mathcal{H}\models_\omega (\mathbf{0}\ra \mathbf{1}).$$
 However, $\mathbf{0}\ra \mathbf{1}$ cannot be proved in the Finitary Injectivity Logic. In fact, the graph
 $$G=\coprod_{n\in\mathbb{N}}C_n$$ demonstrates that $\mathcal{H}\not\models (\mathbf{0}\ra \mathbf{1})$.
 }\es

\end{document}